\documentclass[a4paper,11pt,abstracton,headsepline,bibtotoc]{scrartcl}

\RequirePackage{amsmath}                        
\RequirePackage{amssymb}                        
\RequirePackage{latexsym}                       
\RequirePackage{stmaryrd}                       
\RequirePackage[all]{xy}                        
\UseTips                                        
\newdir^{c}{{}*!/-3pt/\dir^{(}}
\newdir_{c}{{}*!/-3pt/\dir_{(}}

\RequirePackage{xspace}                         
\RequirePackage{calc}                           
\RequirePackage{ifthen}                         

\RequirePackage[dvipsnames,usenames]{color}  
\RequirePackage{hyperref}
\hypersetup{
bookmarks,
bookmarksdepth=3,
bookmarksopen,
bookmarksnumbered,
pdfstartview=FitH,
colorlinks,backref,hyperindex,
linkcolor=RawSienna,
anchorcolor=BurntOrange,
citecolor=OliveGreen,
filecolor=BlueViolet,
menucolor=Yellow,
urlcolor=OliveGreen
}

\usepackage[ansinew]{inputenc}

\KOMAoptions{DIV=13,draft=false,twoside=false,fontsize=11pt,headinclude=false,footinclude=false}%
\def\m?#1{}
    \def\g?#1{}
    \def\mf?#1{}
    \def\forget#1{}

\renewcommand{\sectionautorefname}{Section}

\makeatletter
\@ifdefinable\equationname{\let\equationname\equationautorefname}
\def\equationautorefname~#1\@empty\@empty\null{(#1\@empty\@empty\null)}%
\@ifdefinable\AMSname{\let\AMSname\AMSautorefname}
\def\AMSautorefname~#1\@empty\@empty\null{(#1\@empty\@empty\null)}%
\@ifdefinable\itemname{\let\itemname\itemautorefname}
\def\itemautorefname~#1\@empty\@empty\null{(#1\@empty\@empty\null)%
}%
\makeatother

\makeatletter
\renewcommand{\theenumi}{\alph{enumi}}

\renewcommand{\theenumii}{\roman{enumii}}

\renewcommand{\p@enumii}{\theenumi$\m@th\vert$}
\renewcommand{\theenumiii}{\arabic{enumiii}}

\renewcommand{\p@enumiii}{\theenumi.\theenumii.}

\renewcommand{\labelitemi}{$\m@th\circ$}
\renewcommand{\labelitemii}{$\m@th\diamond$}
\renewcommand{\labelitemiii}{$\m@th\star$}
\renewcommand{\labelitemiv}{$\m@th\cdot$}
\makeatother

\RequirePackage{amsthm}
\RequirePackage{aliascnt}
\newcommand{\basetheorem}[3]{%
    \newtheorem{#1}{#2}[#3]
    \newtheorem*{#1*}{#2}
    \expandafter\def\csname #1autorefname\endcsname{#2}
}%
\newcommand{\maketheorem}[3]{%
    \newaliascnt{#1}{#3}
    \newtheorem{#1}[#1]{#2}
    \aliascntresetthe{#1}
    \expandafter\def\csname #1autorefname\endcsname{#2}
    \newtheorem{#1*}{#2}
}%
\theoremstyle{plain}   

\basetheorem{theorem}{Theorem}{section}

\maketheorem{proposition}{Proposition}{theorem}
\maketheorem{corollary}{Corollary}{theorem}
\maketheorem{lemma}{Lemma}{theorem}
\maketheorem{conjecture}{Conjecture}{theorem}

\theoremstyle{definition}    

\maketheorem{definition}{Definition}{theorem}
\maketheorem{defprop}{Definition-Proposition}{theorem}

\theoremstyle{remark}    

\maketheorem{remark}{Remark}{theorem}
\maketheorem{remarks}{Remarks}{theorem}
\maketheorem{example}{Example}{theorem}
\maketheorem{examples}{Examples}{theorem}
\maketheorem{question}{Question}{theorem}

\numberwithin{equation}{section}

\newcommand{\ie}{\textit{i.e.}}
\newcommand{\cf}{\textit{cf.}}
\newcommand{\blank}{\phantom{m}}

\newcommand{\rank}{\operatorname{rank}}

\newcommand{\dirsum}{\oplus}
\newcommand{\degloc}{\operatorname{degloc}}
\newcommand{\ndsupp}{\operatorname{ndsupp}}
\newcommand{\Fr}{{F}}
\newcommand{\Primefield}{\mathbb{F}}
\newcommand{\CO}{\mathcal{O}}
\newcommand{\CC}{\mathcal{C}}
\newcommand{\CJ}{\mathcal{J}}

\newcommand{\CS}{\mathcal{S}}

\newcommand{\Hom}{\operatorname{Hom}}

\newcommand{\End}{\operatorname{End}}

\newcommand{\tensor}{\otimes}
\newcommand{\Spec}{\operatorname{Spec}}
\newcommand{\defeq}{\stackrel{\scriptscriptstyle \operatorname{def}}{=}}

\newcommand{\Image}{\operatorname{img}}

\renewcommand{\to}[1][]{\xrightarrow{\ #1\ }}

\renewcommand{\frm}{\frak{m}}
\newcommand{\frp}{\frak{p}}
\newcommand{\fra}{\frak{a}}

\renewcommand{\phi}{\varphi}

\newcommand{\frb}{\frak{b}}

\newcommand{\id}{\operatorname{id}}
\newcommand{\supp}{\operatorname{Supp}}
\newcommand{\Ann}{\operatorname{Ann}}
\newcommand{\BF}{\mathbb{F}}

\renewcommand{\th}{\ensuremath{^{\text{th}}}}

\newcommand{\Ca}{{C}}

\reversemarginpar
\makeatletter
\advance\marginparwidth by \ta@bcor
\makeatother
\newcounter{themargin}
\def\m?#1{\textcolor{Mahogany}{\textbf{???$^{\text{\arabic{themargin}}}$}}{\marginpar{\footnotesize\color{Mahogany}\fbox{\parbox{\marginparwidth}{\textbf{mnl --- \arabic{themargin} ---}\addtocounter{themargin}{1}\\ #1}}} \immediate\write16{}%
\immediate\write16{Warning: There was still a question mark . . . }%
\immediate\write16{}}}
\def\mf?#1{\textcolor{Mahogany}{\footnotesize\newline{\color{Mahogany}\fbox{\parbox{\textwidth-5pt}{\textbf{mnl: } #1}}}}}
\definecolor{WarmDarkGray}{rgb}{0.25,0.2,0.18}

\begin{document}

\title{Test ideals via algebras of $p^{-e}$-linear maps}
\author{Manuel Blickle}
\maketitle
\footnotetext{to appear in Journal of Algebraic Geometry}
\footnotetext{\textbf{Address:} Johannes Gutenberg-Universität Mainz, FB08, Institut für Mathematik, 55000 Mainz, Germany}
\footnotetext{\textbf{email:} \texttt{manuel.blickle@gmail.com}}
\footnotetext{\textbf{keywords:} Frobenius, $p^{-1}$-linear maps, Cartier operator, tight closure, test ideals}
\footnotetext{\textit{2010 Mathematics subject classification.} Primary: 13A35; Secondary: 14G17}

\begin{abstract}
    Building on  ideas of a recent preprint of Schwede \cite{schwede_test_2009}, combined with the theory developed by Böckle and the author in \cite{BliBoe.CartierFiniteness}, we study test ideals by viewing them as minimal objects in a certain class of modules, called $F$-pure modules, over algebras of $p^{-e}$-linear operators. We develop the basics of a theory of $F$-pure modules and show an important structural result, namely that $F$-pure modules have finite length. This result is then linked to the existence of test ideals and  leads to a simplified and generalized treatment, also allowing us to define test ideals in non-reduced settings.

    Combining our approach with an observation of Anderson in \cite{AndersonL} on the contracting property of $p^{-e}$-linear operators yields an elementary approach to test ideals in the case of affine $k$-algebras, where $k$ is an $F$-finite field. As a byproduct one obtains a short and completely elementary proof of the discreteness of the jumping numbers of test ideals in a generality that extends most cases known so far, \cf~\cite{BliSchTakZha_DisccRat}, in particular one obtains results beyond the $\mathbb{Q}$-Gorenstein case.
\end{abstract}

\section{Introduction}
Throughout this article we fix a noetherian ring $R$ of positive characteristic $p$. We often assume that the Frobenius, \ie~the $p\th$ power map, is a finite map. Test ideals play an important role in the classification of singularities in positive characteristic. They are analogues of the multiplier ideals which are a fundamental tool in the birational geometry of varieties in characteristic zero. However, test ideals are in many respects more subtle as they also encode arithmetic properties of the variety. Test ideals developed out of the theory of tight closure \cite{HH90}, their connection to multiplier ideals was worked out in \cite{Smith.multipl,Hara.GeomTight} and extended to the relative setting in \cite{HaraYosh.GenTightMult,Takagi.MultiplierTight}.

Recently, Schwede generalized in \cite{schwede_test_2009} all previous approaches to test ideals by further emphasizing the role that $p^{-e}$-linear maps play in their construction. A $p^{-e}$-linear map on $R$ is an additive map $\phi: R \to R$ that satisfies $\phi(r^{p^e}s)=r\phi(s)$ for all $r,s \in R$. The connection of $p^{-e}$-linear maps with test ideals has implicitly been used before, but Schwede in \cite{schwede_test_2009} and before in \cite{schwede_centers_2008,schwede_f-adjunction_2009} successfully emphasizes this viewpoint. His idea is to associate a test ideal to a pair $(R,\CC)$ where $R$ is a noetherian, reduced $F$-finite ring, and $\CC$ is a \emph{Cartier algebra over $R$}, that is an algebra of $p^{-e}$-linear operators on $R$, \cf~\autoref{d.RCartAlgebra}. Now, roughly  speaking (and slightly incorrectly, see \autoref{d.test}), the test ideal $\tau(R,\CC)$ associated to the pair $(R,\CC)$ is the smallest ideal of $R$ that is also a $\CC$-submodule of $R$ and and which agrees with $R$ generically. By choosing appropriate algebras of $p^{-e}$-linear operators this definition recovers all previously considered notions of test ideals for pairs $(R,\fra^t)$ or triples $(R,\Delta,\fra^t)$ where $\fra$ is an ideal, $t \in \mathbb{R}_{>0}$, and $\Delta$ is an effective $\mathbb{Q}$-divisor, see~\autoref{s.CartierAlgebrasExamples}.

The aim of this manuscript is to use the techniques of \cite{BliBoe.CartierFiniteness} to generalize this viewpoint on test ideals further by considering pairs $(M,\CC)$ where $M$ is a finitely generated $R$-module, $R$ is an arbitrary Noetherian and $F$-finite ring (not necessarily reduced), and $\CC=\oplus_{e=0}^\infty \CC_e$ is an $\mathbb{N}$-graded algebra of $p^{-e}$-linear operators acting on $M$. Our aim is to construct test modules in this setting as minimal objects in a suitable category of $\CC$-submodules of $M$. This category is called $F$-pure $\CC$-submodules and consist of those submodules $N \subseteq M$ such that $\CC_+N=N$ where $\CC_+=\oplus_{e=1}^\infty \CC_e$. The generality of this approach leads to a simplified theory which also includes a candidate for the test ideal in the non-reduced setting. We generalize the notions of (strong) $F$-regularity, and in particular, $F$-purity to this setting. In the special case that the algebra $\CC$ is principally generated, the main result in \cite[Theorem 4.6]{BliBoe.CartierFiniteness} states that in the category of coherent $F$-pure $\CC$-modules all objects have finite length. Here we show that this finite length condition is for general $\CC$ intimately tied to the existence of test modules:
\begin{theorem*}[\cf~\autoref{t.finitelength}]
Let $R$ be a Noetherian ring and $\CC$ an $R$-Cartier algebra, $M$ a finitely generated $\CC$-module. Then all chains of $F$-pure $\CC$-submodules of $M$ are finite if and only if the test modules $\tau(N,\CC)$ exist for all $\CC$-submodules $N \subseteq M$.
\end{theorem*}
As the statement of this result shows, the crucial point is the existence of the test modules as defined above. This is established not in complete generality, but for many important cases. After we explained the basic notions and results of our theory of $F$-pure modules in \autoref{s.Cart} and \autoref{s.TestModules} we show in \autoref{s.existenceTest} the following result
\begin{theorem*}[\cf~\autoref{t.ExistenceTestModule}, \autoref{t.existenceTestModule2}]
Let $R$ be a $F$-finite Noetherian ring, and $\CC$ be an $R$-Cartier-algebra and $M$ a coherent left $\CC$-module. In each of the following cases does the test module $\tau(M,\CC)$ exist.
\begin{enumerate}
\item $M \subseteq R$
\item $\CC$ is principally generated
\item $R$ is of finite type over an $F$-finite field
\end{enumerate}
\end{theorem*}
The proof in the first case is quite similar to the classical proofs of existence of test elements. The second item follows from the theory developed in \cite{BliBoe.CartierFiniteness}. The final part is new and a byproduct of another viewpoint on the study of test ideals which uses the contracting property of $p^{-e}$-linear maps, as nicely explained by Anderson in \cite{AndersonL}. Emphasizing this contracting property of $p^{-e}$-linear maps we obtain a completely elementary proof of the existence of test modules as well as the discreteness of jumping numbers for test modules of triples $\tau(M,\CC,\fra^t)$ with $\fra$ an ideal of $R$ and $t \in \mathbb{R}_{\geq 0}$ in the case that $R$ is of finite type over an $F$-finite field, and $\CC$ satisfies a certain boundedness condition, \cf~\autoref{s.elementary} for definitions and details:
\begin{theorem*}[\cf~\autoref{t.TestExistsAffineKalg}, \autoref{t.GaugeBoundDiscrete}]
Let $R$ be of finite type over an $F$-finite field, $\fra \subseteq R$ and ideal, $\CC$ an algebra of $p^{-e}$-linear maps and $M$ a $R$-finitely generated $\CC$-module.
\begin{enumerate}
\item The test module $\tau(M,\CC,\fra^t)$ exists for all $t \geq 0$.
\item If $\CC$ is a finitely generated $R$-algebra, then there are only finitely many submodules in the collection $\tau(M,\CC,\fra^t)$ for $0 \leq t \leq T$ for any fixed $T >0$.
\end{enumerate}
\end{theorem*}
The second statement is our analog of the discreteness of $F$-jumping numbers. It greatly extends most cases obtained before in \cite{BliMusSmi.DisRat,BliSchTakZha_DisccRat}. However, the proof we present here is completely elementary as it avoids any use of hard machinery such as Schwede's $F$-adjunction that was a crucial ingredient in \cite{BliSchTakZha_DisccRat}, or the $D$-module techniques of \cite{BliMusSmi.HypF}. We also note that even in characteristic zero there is little known on the discreteness of the jumping numbers of the multiplier ideals beyond the $\mathbb{Q}$-Gorenstein case (see however \cite{urbinati_discrepancies_2010}).

There are many directions for new research that emerge from the theory we introduce here, most prominently would be a detailed study of the category of $F$-pure modules associated to any pair $(R,\CC)$. However, in order to keep this paper at a reasonable size we decided to only give the basic constructions and ideas necessary to set up the theory and quickly show the results mentioned above. This is done at the cost of omitting many other interesting aspects. For example, we did not touch on generalizations of classical properties of test ideals, such as subadditivity, nor do we study the new invariants our theory produces (\textit{e.g.}~the length in the category of $F$-pure modules) and their connection to previously studies invariants of singularities. By including a list of questions in the final \autoref{s.question} we hope to stimulate further research in these directions.

\subsection*{Acknowledgements}
I thank Karl Schwede for stimulating discussions and for sending me an early draft of his manuscript for \cite{schwede_test_2009} which initiated this paper. I also thank Daisuke Hirose for pointing out a small mistake in an earlier version of this paper. I also thank the anonymous referee for a very careful reading and many valuable suggestions.

The research for this paper was conducted while I was supported by a Heisenberg Fellowship of the DFG and by the SFB/TRR45.

\subsection*{Notation and Conventions}
We fix a prime number $p$. Let $\Primefield=\Primefield_p$ be the finite field with $p$ elements. We will work in the affine setting for convenience, however, the theory we introduce does localize well, such that the constructions and results are valid for general noetherian schemes. All rings $R$ will be noetherian $\Primefield_p$-algebras and we denote by $\Fr$ the $p$\th-power Frobenius, \ie~the $p$\th-power map on $R$. For the most part we also assume that $R$ is $\Fr$-finite, i.e. the Frobenius $\Fr$ is a finite map of schemes. For an $R$-module $M$ we denote by $\Fr_*M$ the additive group $M$ with the $R$ structure given by $r\cdot m = r^{p}m$.

\section{Algebras of $p^{-e}$-linear maps}
\label{s.Cart}
The basic objects we are concerned with are $p^{-e}$-linear maps on $R$-modules. This topic has a long history and the prime example is the Cartier operator on the dualizing module of $R=k[x_1,\ldots,x_n]$ with $k$ a perfect field, \cf~\cite{Cartier57}. Another typical example of a $p^{-1}$-linear map is a splitting of the $p$-linear Frobenius map, see \cite{BrionKumar.FrobSplit}.
\begin{definition}
\label{d.CartLinearMap}
A $p^{-e}$-linear map of $R$-modules is an additive map $\phi_e : M \to N$ that satisfies $\phi_e(r^{p^e}m)=r\phi_e(m)$ for all $r \in R$ and $m \in M$; equivalently, it is an $R$-linear map $\phi_e: \Fr^e_*M \to N$. We denote by $\Hom_e(M,N)=\Hom_R(\Fr^e_*M,N)$ the abelian group of $p^{-e}$-linear maps.
\end{definition}
Even though a $p^{-e}$-linear map $\phi_e: M \to N$ is not $R$-linear, its image $\phi_e(M)$ is an $R$-submodule of $N$, since $r\phi_e(m)=\phi_e(r^{p^e}m)$. However, the kernel is only a $R$-submodule of $\Fr_*^eM$ (equiv. an $R^{p^e}$-submodule of $M$). One can compose a $p^{-e}$-linear map $\phi_e$ and a $p^{-e'}$-linear map $\phi_{e'}$ in the obvious way as additive maps, the result is a $p^{-(e+e')}$-linear map. Viewing $\phi_e$ as an $R$-linear map $\Fr^e_*M \to N$ and $\phi_{e'}: \Fr^{e'}_*N \to L$, this is just the composition $\phi_{e'} \circ \Fr^{e'}_*\phi_{e}  :\Fr^{*(e+e')}M \to L$.
\begin{definition}
\label{d.RCartAlgebra}
Let $R$ be a commutative noetherian ring of characteristic $p>0$. An \emph{$R$-Cartier-algebra} (or $R$-algebra of Cartier type) is an $\mathbb{N}$-graded $R$-algebra $\CC := \oplus_{e \geq 0} \CC_e$ such that for each $r \in R$ and $\phi_e \in \CC_e$ one has $r \cdot \phi_e = \phi_e \cdot r^{p^e}$. The $R$-algebra structure of $\CC$ is nothing but a ring homomorphism $R \to \CC$, in fact, the image of $R$ lies in $\CC_0$ since $1 \in R$ maps to $1 \in \CC_0$. We will assume that this structural map $R \to \CC_0$ is surjective.
\forget{\footnote{It seems convenient to assume here that $R$ surjects onto $\CC_0$, but I am not sure if it is absolutely necessary. Another generalization would be to assume that the Frobenius $F$ acts on the $R$-algebra $\CC_0$ and that one has for every $\phi \in \CC_e$ and $s \in \CC_0$ that $s \cdot \phi=\phi \cdot F^e(s)$. One other interesting point of view may be to consider, for $\CC$ an $R$-Cartier algebra, the (graded) tensor product $R[F] \tensor_R \CC$, where we should think of $F$ sitting in degree $-1$. If one considers this, it would probably be better to let $F^e$, which is a $p^e$-linear map, have degree 1, and a $p^{-e}$ linear map to have degree $-e$.}}
\end{definition}
Note that the algebra $\CC$ is generally not commutative, not even $\CC_0$ or $R$ are central. Each graded piece $\CC_e$ of $\CC$ is naturally a $R-R$-bi-module, and the defining relation shows that the left structure is determined by the right $R$-structure.
\begin{example}
Let $M$ be an $R$-module, then we denote by
\[
    \End_*(M)=\oplus_{e \geq 0} \End_e(M)
\]
the \emph{algebra of Cartier linear operators on $M$}. According to our definition, this is generally $\emph{not}$ an $R$-Cartier-algebra, since $\End_0(M)=\End_R(M)$ is generally too big. However, specializing to the case of $M=R$ we define
\[
    \CC_R := \End_*(R)
\]
the \emph{algebra of Cartier linear operators on $R$}, which is an $R$-Cartier-algebra since $\End_R(R)=R$. Again $\CC_{R,0}=R$ is not central, but one has the relation $r\phi_e = \phi_er^{p^e}$ for all $\phi_e \in \CC_{R,e}=\End_e(R)$ and $r \in R$.
\end{example}
\begin{remark}
There is a natural map
\[
    \CC_R \to \End_\mathbb{Z}(R)
\]
given by forgetting the $p^{-e}$-linearity. In reasonable circumstances this map is injective.\footnote{Similar to \cite{EmKis.Fcrys} Lemma 1.4.1 one may show that if all minimal primes of $R$ have height $\geq 1$, then the natural map $\CC_R \to \End_\mathbb{Z}(R)$ is injective}
\forget{The following lemma clarifies when this map in an injection.
\begin{lemma}[cf.~\protect{\cite{EmKis.Fcrys}} Lemma 1.4.1]
If all minimal primes of $R$ have height $\geq 1$, then the natural map $\CC_R \to \End_\mathbb{Z}(R)$ is injective.
\end{lemma}
\begin{proof}
The condition ensures that each maximal ideal $\frm$ of $R$ contains a non-zero-divisor $r$. Suppose $\sum_{i=0}^e \phi_i \in \CC_R$ with $\phi_i \in \CC_{R,i}$ is an element the kernel. This means, for all $s \in R$ and we have $\sum_{i=1}^e \phi_i(s)=0$. We have to show that each $\phi_i$ is zero. Since a map $\phi_i: \Fr_*R \to R$ is zero if and only if its image, which is an ideal in $R$, is zero, we may check this locally, and hence assume that $R$ is local and $r$ is a non-zero divisor on $R$. In fact, by assuming that the degree (in $e$) of $\sum_{i=0}^e \phi_i \in \CC_R$ is minimal amongst all elements in the kernel, it is enough to show that $\phi_e=0$. Replacing $s$ by $r^{mp^e}s$ and using the $p^{-i}$-linearity of each $\phi_i$ we get $\sum_{i=0}^e r^{mp^{e-i}}\phi_i(s)=0$ for all $s \in R$ and $m > 0$. Dividing by $r$ we get
\[
    0 = \sum_{i=0}^{e} r^{(p^{e-i}-1)}\phi_i(s)\, .
\]
In this sum, all but the last term are divisible by $r^{m(p-1)}$. Hence the last term $\phi_e(s)$ has to be divisible by $r^{m(p-1)}$ as well. Since $m > 0$ was arbitrary and $r$ is a non-zero-divisor it follows that $\phi_e(s)=0$ for all $s$, hence $\phi_e=0$.
\end{proof}
}
An example where it is \emph{not} injective, is $R= \Primefield_p$, where we have $\End_{\Primefield_p}(\Primefield_p)=\Primefield_p$ but $\CC_{\Primefield_p} \cong \dirsum_{e \geq 0} \Primefield_p$.
\end{remark}

\begin{remark}
For clarification, note that $\End_e(R) = \Hom_{R}(\Fr^e_*R,R)$ is a $R$-$R$-bi-module. The right and left $R$-structures are given by pre- and post-composition of $p^{-e}$-linear maps with multiplication by $r \in R$. That means for $r \in R$, and $\phi \in \End_e(R)$ we have two maps $\phi \cdot r$ and $r \cdot \phi$ in $\End_e(R)$ defined for $s \in R$: $(\phi \cdot r)(s)=\phi(rs)$ and $(r \cdot \phi)(s)=r\phi(s)=\phi(r^{p^e}s)$. The right $R$-structure might also be thought of as the $\Fr^e_*R$-module structure on the first entry, after one identifies $R$ and $\Fr^e_* R$ as rings. For $\Hom_{R}(\Fr^e_*R,R)$ with this right $R$-module structure one traditionally writes $\Fr^{e\flat}R$. Hence $\CC_R \cong \oplus_{e \geq 0} \Fr^{e\flat}R$.
\end{remark}
In this article we study left modules over some $R$-Cartier-algebra $\CC$, with a particular view towards the theory of test ideals, and the resulting generalizations. Before doing so, we discuss a little bit on the structure of some $R$-Cartier-algebras we will be considering. We begin with the principally generated case, which was studied in detail in \cite{BliBoe.CartierFiniteness}. It is the prevalent case if the underlying ring $R$ is Gorenstein.
\begin{example}
Let $R$ be a noetherian ring. Let $R[\Fr]$ be the non-commutative ring arising from $R$ by adjoining a new variable $\Fr$ subject to the relations $r^p\Fr-\Fr r$ for all $r \in R$. As a left $R$-module one has an isomorphism $R[\Fr] \cong \oplus_{e \geq 0} R\Fr^e$. Right from the definition one checks that the opposite algebra of $R[\Fr]$, defined as $R\langle \Ca \rangle \defeq R[\Fr]^{op}$ is an $R$-Cartier-algebra. As \emph{right} $R$-modules one has $R\langle \Ca \rangle = \oplus_{e \geq 0} \Ca^eR$. This algebra is the prototype of an algebra of Cartier linear operators. Of course, it may also be obtained by adjoining to $R$ the variable $\Ca$ subject to the relation $r\Ca-\Ca r^p$ for all $r \in R$. The theory of Cartier modules as developed in \cite{BliBoe.CartierCrys,BliBoe.CartierFiniteness} is, in fact, a theory of left modules over $R\langle \Ca \rangle$.
\end{example}
Next we discuss the structure of the $R$-Cartier-algebra $\CC_R$ in some detail. In particular we are interested in its relation to the algebra $R\langle \Ca \rangle$ or its Veronese subalgebras $R\langle \Ca^n \rangle$ for some $n \geq 1$.
\begin{example}
Let $R$ be Gorenstein and $F$-finite, then
    \[
        \CC_{R,e}=\Hom_{R}(\Fr^e_*R,R) \cong \Fr^{e\flat}R \cong \omega^{1-p^e}_R
    \]
is locally free of rank one as a right $R$-module. Denote by $\phi_e$ a (local) right $R$-generator of $\CC_{R,e}$. Then $\phi_e^{e'}$ is a right $R$-generator of $\CC_{R,ee'}=\Hom_{R}(\Fr^{ee'}_*R,R) \cong \omega^{1-p^{ee'}}_R$. Hence $\CC_R$ is principally generated by $\phi_1$ as an $R$-algebra. Mapping $\Ca$ from the preceding example to the generator $\phi_1$, yields an isomorphism of $R$-algebras $R\langle C\rangle \to \CC_R$.
\end{example}
What this example shows is that if $R$ is reduced, then as a (right) $R$-algebra we have that $\CC_R$ is generically principally generated, namely on any open set $U$ where $X$ is Gorenstein.
More generally, if $R$ is $F$-finite (and normal), then by \cite{Gabber.tStruc} it has a dualizing complex, and hence by \cite{MehtaRam} one has (at least Zariski locally) $\Hom_{R}(\Fr^e_*R,R) \cong \Fr^e_*R((1-p^e)K_X)$ as $\Fr^e_*R$-modules, so $\CC_{R,e} \cong \omega_R^{1-{p^e}}$ as right $R$-modules. In general, if $R$ is not Gorenstein, this shows that $\CC_R$ is not principally generated. In fact, it may not even be finitely generated as the following explains.

Let $R$ be complete, local, and $F$-finite. Then $\CC_R$ is isomorphic to  the opposite of the algebra of Frobenius-linear operators on the injective hull of the residue field of $k$, \cf~\cite{BliBoe.CartierFiniteness}. An example of Katzman \cite{Katzman.NonFG} shows that the latter need not be finitely generated.

Note that there is always a canonical map from $R[\Fr]$ to the algebra of Frobenius linear operators $\End^*(R)=\oplus_{e\geq 0}\Hom(R,\Fr^e_*R)$ sending $\Fr$ to the Frobenius. Unless $R$ is zero dimensional and of finite type over $\Primefield_p$ this map is an injection \cite{EmKis.Fcrys}. In the dual situation of $p^{-e}$-linear operators, there is generally no canonical map $R\langle \Ca \rangle \to \CC_R$. However, in many cases (for example $R$ is $F$-finite and sufficiently affine), there is a (somewhat) canonical map $R\langle \Ca \rangle \to \End_*(\omega_R)$ given by sending $\Ca$ to the dual (under Grothendieck-Serre duality) of the Frobenius on $R$.

\subsection{Nilpotence}
In the case of a single $p^{-e}$-linear operator, the notion of nilpotence was treated thoroughly in \cite{BliBoe.CartierCrys,BliBoe.CartierFiniteness}. We slightly generalize here some basic constructions of $\cite{BliBoe.CartierFiniteness}$ from the case of a single operator to that of an algebra of operators. For this let $\CC$ be an $R$-Cartier-algebra. Denote by $\CC_+=\oplus_{e \geq 1} \CC_e$ the positively graded part of $\CC$. This is a two-sided ideal in $\CC$ and we denote by
\[
    \CC_+^n=\{ \phi_1\cdot\ldots\cdot\phi_n | \phi_i \in \CC_+ \}
\]
its $n$\th~power. For any $\CC$-module $M$ we denote by $\CC_+^nM$ the $\CC$-submodule of $M$ generated by all $\phi \cdot m$ for $\phi \in \CC_+^n$ and $m \in M$. Since the image of an $R$-module under any $p^{-e}$-linear map is again an $R$-module we have $\CC_+^n M = \sum_{\phi \in \CC_+^n} \phi(M)$.
\begin{definition}
\label{d.nil}
A $\CC$-module $M$ is called \emph{nilpotent} (or $\CC$-nilpotent if $\CC$ is not clear from the context) if $\CC_+^nM=0$ for some $n \geq 0$.
\end{definition}\begin{definition}
\label{d.coh}
A $\CC$-module $M$ is called \emph{coherent} if it is finitely generated as an $R$-module.
\end{definition}
\begin{lemma}
The subcategory of $\CC$-nilpotent $\CC$-modules, as well as the subcategory of coherent $\CC$-modules are Serre subcategories of the category of $\CC$-modules, \ie~both are abelian subcategories which are closed under extensions.
\end{lemma}
\begin{proof}
All statements are clear except that the extension of nilpotent $\CC$-modules is again nilpotent. This follows as in \cite[Lemma 2.11]{BliBoe.CartierFiniteness}.
\end{proof}
Before proceeding we explain how $\CC$-modules behave under localization. Let $S$ be a multiplicative subset of $R$. Each $\phi_e \in \CC_e$ naturally acts on $S^{-1}M$ by the formula $\phi_e(m/s) = \phi_e(ms^{p^e-1})/s$. Hence $S^{-1}M$ naturally has the structure of a left $S^{-1}\CC$-module, where $S^{-1}\CC=\CC \tensor_R S^{-1}R$. The latter may also be viewed as the localization of the non-commutative $R$-algebra $\CC$ at the image of the multiplicative set $S$ in $\CC$, via the $R$-algebra structure. Even though this multiplicative set is generally not central in $\CC$ the localization works here just as in the commutative case since one has the relation $S^{-1}\CC_e=\CC_eS^{-p^e}=\CC_eS^{-1}$. Summarizing we obtain:
\begin{lemma}
\label{t.CC+CommutesLocalization}
Let $\CC$ be an $R$-Cartier-algebra, $M$ a left $\CC$-module, and $S \subseteq R$ a multiplicative subset of $R$. Then $S^{-1}M$ is naturally a left module over the $S^{-1}R$-Cartier algebra $S^{-1}\CC$. Furthermore one has
\[
    S^{-1}(\CC_+M) = (S^{-1}\CC)_+(S^{-1}M) \\
\]
\end{lemma}
\begin{proof}
    Just compute:
    \[
        S^{-1}(\CC_+M)= \sum_{e \geq 1} S^{-1}(\CC_e M)=\sum_{e \geq 1} \CC_e(S^{-p^e} M)= \sum_{e \geq 1} \CC_e(S^{-1}M) = (S^{-1}\CC)_+(S^{-1}M)
    \]
\end{proof}
This preceding lemma ensures that the key constructions that are to follow behave well with respect to localization.
\begin{lemma}
Let $M$ be a coherent $\CC$-module. There is a unique submodule $M_{nil}$ such that
\begin{enumerate}
\item $M_{nil}$ is nilpotent, and
\item $\overline{M}=M/M_{nil}$ does not have nilpotent $\CC$-module quotients.
\end{enumerate}
\end{lemma}
\begin{proof}
Let $M_e$ be the largest $\CC$-submodule of $M$ such that $\CC^e_+ M_e=0$. This exists since the sum of two such also has this property. Since $M$ is coherent, the union $M_{nil} := \bigcup_e M_e$ stabilizes after finitely many steps. It follows that $M_{nil}$ is in fact $\CC_+$-nilpotent. The rest of the proof is analogous to the one in \cite[Lemma 2.12]{BliBoe.CartierFiniteness}.
\end{proof}
The following Proposition is absolutely critical to our theory. It may be viewed as a vast generalization of a result on the nilpotence of locally nilpotent Artinian modules with a left Frobenius action of Hartshorne and Speiser \cite{HaSp}. The proof that follows is an adaption of the one in \cite[Lemma 13.1]{Gabber.tStruc} in the case of a single operator.
\begin{proposition}
Let $M$ be a coherent $\CC$-module. The descending chain
\[
    M \supseteq \CC_+(M) \supseteq \CC_+^2(M) \supseteq \CC_+^3(M) \supseteq \ldots
\]
of $\CC$-submodules stabilizes.
\end{proposition}
\begin{proof}
Each $\CC_+^e(M)$ is a coherent $\CC$-submodule. By \autoref{t.CC+CommutesLocalization} we have for each multiplicative set $S$ that $S^{-1}(\CC_+^e(M))=(S^{-1}\CC)_+^e(S^{-1}M)$. Define
\[
    Y_e := \supp (\CC^e_+(M)/\CC_+(\CC^e_+(M)))
\]
which form a descending sequence of closed subsets of $X=\Spec R$. By noetherian-ness this descending sequence must stabilize. After truncating we may assume that for all $e$ we have $Y:=Y_e=Y_{e+1}$. The statement that the above chain stabilizes means precisely that $Y$ is empty. Let us assume otherwise and let $\frm$ be the generic point of a component of $Y$. Localizing at $\frm$ we may assume that $(R,\frm)$ is local and that $\{\frm\}=Y=\supp (\CC^e_+(M)/\CC_+(\CC^e_+(M)))$ for all $e$. In particular, for $e=0$ we get that there is a $k$ such that $\frm^k M \subseteq  \CC_+(M)$. Then, for any $x \in \frm^k$
\[
    x^2M \subseteq x\frm^k M \subseteq x\CC_+M \subseteq \sum_{e \geq 1} \CC_e(x^{p^e}M)\subseteq \CC_+(x^2M)
\]
and iterating $x^2M \subseteq \CC_+^e(M)$ for all $e$. Hence $\frm^{k(b-1)}M \subseteq \CC_+^eM$ for all $e$ where $b$ is the number of generators of $\frm^k$. Hence the original chain stabilizes if and only if the chain $\CC_+^eM/\frm^{k(b-1)}M$ does. But the latter is a chain in the finite length module $M/\frm^{k(b-1)}M$.
\end{proof}
\begin{corollary}
\label{t.Underline}
Let $M$ be a coherent $\CC$-module. There is a unique $\CC$-submodule $\underline{M}$ such that
\begin{enumerate}
\item the quotient $M/\underline{M}$ is nilpotent, and
\item $\CC_+\underline{M}=\underline{M}$ (\ie~$\underline{M}$ does not have nilpotent quotients).
\end{enumerate}
\end{corollary}
\begin{proof}
The stable member $\underline{M}:= \CC_+^e(M)$ for $e \gg 0$ has all the desired properties.
\end{proof}
Note that even though this is hidden by the notation, the functors sending $M$ to $M_{nil}$, $\underline{M}$, and $\overline{M}$ all depend on the algebra $\CC$. Hence, if confusion may occur we will denote these by $M_{\CC-nil}$, $\underline{M}_{\CC}$, and $\overline{M}^{\CC}$. This innocent looking operation of passing from $M$ to $\underline{M}$ is absolutely crucial in our treatment. The relative technical simplicity, and wide scope of our approach to test ideals which follows below owes much to a systematic use of this operation. We make the following definition.
\begin{definition}
Let $\CC$ be an $R$-Cartier-algebra and let $M$ a left $\CC$-module. The pair $(M,\CC)$ (or just $M$, if $\CC$ is clear from the context) is called \emph{$F$-pure}, if $\CC_+M=M$.
\end{definition}
This terminology may be justified by observing that if $M=R$, then $(R,\CC)$ is $F$-pure if and only if $\CC$ contains a splitting of some power of the Frobenius on $R$, \cf~\autoref{t.FpureCharacterization} below; in particular, the Frobenius map $F$ on $R$ is a pure map of rings. Some simple reformulations of $F$-purity are:
\begin{lemma}
Let $\CC$ be an $R$-Cartier algebra and $M$ a left coherent $\CC$-module. Then the following conditions are equivalent:
\begin{enumerate}
\item $(M,\CC)$ is $F$-pure.
\item $\CC^e_+M=M$ for all $e \geq 1$.
\item $\underline{M}_{\CC}=M$
\item $M$ does not have $\CC$-nilpotent quotients.
\end{enumerate}
\end{lemma}

The following important corollary says that, up to nilpotence, we may replace any $\CC$-module by its (unique) maximal $F$-pure submodule.
\begin{corollary}
Let $M$ be a coherent $\CC$-module. Then $\underline{M}_{\CC}$ is the largest $\CC$-submodule $N$ of $M$ such that $\CC_+N=N$, \ie~it is the largest $F$-pure submodule. Furthermore, $M/\underline{M}_{\CC}$ is the maximal nilpotent quotient of $M$.
\end{corollary}
\begin{proof}
Since the sum of $\CC$-submodules $N_i$ of $M$ with the property $\CC_+N_i=N_i$ also has this property it is clear that there is a (unique) largest $\CC$-submodule $N$ with $\CC_+N=N$. Since $N \subseteq M$ it follows that $N=\CC_+N \subseteq \CC_+M$, and iterating we get $N \subseteq \CC_+^e M$ for all $e \geq 0$. Hence $N \subseteq \underline{M}_{\CC}$ and $\underline{M}_{\CC}$ is indeed the largest such submodule. The final statement follows immediately from the construction of $\underline{M}_{\CC}$ in \autoref{t.Underline}.
\end{proof}
As an immediate consequence of \autoref{t.CC+CommutesLocalization} we get.
\begin{lemma}
\label{t.CommutesLocalization}
Let $\CC$ be an $R$-Cartier-algebra, $M$ a left $\CC$-module, and $S \subseteq R$ a multiplicative subset of $R$. Then
\[
    S^{-1}\underline{M}_{\CC} = \underline{S^{-1}M}_{S^{-1}\CC},\blank S^{-1}\overline{M}^{\CC} = \overline{S^{-1}M}^{S^{-1}\CC},\text{ and } S^{-1}({M}_{\CC-nil}) = (S^{-1}M)_{S^{-1}\CC-nil}
\]
\end{lemma}
\begin{proof}
Easy.
\end{proof}

\subsection{Support and restriction to a closed set}
If $M$ is a coherent $\CC$-module, then its support $\supp M = V(\Ann_R M)$ is a closed subset of $\Spec R$. A basic, however important, observation is that if $M$ is $F$-pure, \ie~$\CC_+M=M$, then the support is reduced, \ie~the annihilator of $M$ is a radical ideal.
\begin{lemma}
\label{t.SupportReduced}
Let $M$ be a coherent $F$-pure $\CC$-module, then $\supp M$ is reduced, equivalently, $\Ann_R M$ is a radical ideal. In particular, if $\supp M \subseteq \Spec R/I \subseteq \Spec R$ for some ideal $I \subseteq R$, then $IM=0$.
\end{lemma}
\begin{proof}
Let $r^n \in \Ann_R M$. Then, for $f$ with $p^f \geq n$ we have for all $e \geq f$ that $r^{p^e}M=0$. The assumption $\CC_+M=M$ implies $M=\CC^e_+M =\sum_{e \geq f} \CC_{e}M$, and hence
\[
    0=\sum_{e \geq f} \CC_e r^{p^e}M = \sum_{e \geq f} r \CC_eM=rM.
\]
Hence $\Ann_R M$ is a radical ideal.
\end{proof}

If $I \subseteq R$ is an ideal, then the left ideal $\CC  \cdot I$ of $\CC$ is in fact a two-sided ideal, since one has for all $r \in R$ the relation $r \CC_e = \CC_e r^{p^e}$. Hence the quotient $\CC/\CC I \cong \CC \tensor_R R/I$ is again an $R$-Cartier-algebra, indeed it is an $R/I$-Cartier algebra. Clearly, any $\CC$-module $M$ with $IM=0$ may be viewed as a $\CC/\CC I$-module, and vice versa.
\begin{lemma}
\label{t.restrictionClosed}
Let $I$ be an ideal of $R$, and $M$ a left module over the $R$-Cartier-algebra $\CC$. Then the $\CC$-submodule
\[
    \Ann_M I = \{ m \in M | Im=0 \} = \Hom_R(R/I,M)
\]
is naturally a left module over the $R/I$-Cartier algebra $\CC/\CC I$.
\end{lemma}
\begin{proof}
If $Im=0$ then for every $i \in I$ and $\phi \in \CC_e$  we have $i\phi(m)=\phi(i^{p^e}m)=0$. Hence $\Ann_M I$ is a $\CC$-submodule of $M$, and the statement follows from the discussion preceding the lemma.
\end{proof}
Hence, if the support of a coherent $F$-pure $\CC$-module $M$ is contained in $\Spec R/I$ for some ideal $I \subseteq R$, then $\Ann_M I=M$. Consequently, a coherent $F$-pure submodule with support in $\Spec R/I$ may be viewed as a $\CC/\CC I$-module, and vice versa. Of course, the pair $(M,\CC/\CC I)$ is also $F$-pure. Summarizing, we obtain:
\begin{proposition}
Let $I$ be an ideal of $R$, and $\CC$ an $R$-Cartier algebra. Then the category of (coherent) $F$-pure $\CC$-modules which are supported on $\Spec R/I$ is equivalent to the category of (coherent) $F$-pure $\CC/\CC I$-modules.
\end{proposition}

\section{Test modules}
\label{s.TestModules}
In this section we generalize slightly the \emph{ultimate generalization of test ideals} of Schwede in \cite{schwede_test_2009} to the following setting. As generally assumed, $R$ is a noetherian ring over $\Primefield_p$. We will suggest here a definition of test modules\footnote{According to accepted terminology, these should be called \emph{big} test modules, to ease notation in this article I just call them test modules. Note that the classical notions of big test ideals and test ideals agree in many cases, for example if $R$ is $F$-finite and $\mathbb{Q}$-Gorenstein} associated to pairs $(M,\CC)$ consisting of an $R$-Cartier-algebra $\CC$, and a coherent left $\CC$-module $M$. In \cite{schwede_test_2009} the case that $\CC \subseteq \CC_R$ and $M=R$ was considered. The additional freedom we get by the added generality allows us to define test ideals also for non-reduced rings, in the degenerate case (\cf~\cite[Definition 3.8]{schwede_test_2009} for degeneracy), and leads to a simplified and more conceptual theory. The following is the defining property of the test module.
\begin{definition}
\label{d.test}
The \emph{test module} of a pair $(M,\CC)$ is the smallest $\CC$-submodule of $M$ that \emph{agrees generically} with $\underline{M}_{\CC}$, \ie~the smallest $\CC$-submodule $N \subseteq M$ such that for each generic point $\eta$ of a component of $\supp \underline{M}_\CC$ one has $(\underline{M}_{\CC})_\eta = N_\eta$. The test module is denoted by $\tau(M,\CC)$ and it is necessarily contained in $\underline{M}_{\CC}$. If $M=R$ we call $\tau(R,\CC)$ also the \emph{test ideal}.
\end{definition}
The issue is of course whether the test module exists. We will show the existence of test modules assuming that $R$ is $F$-finite, $M$ is finitely generated, and one of the following additional conditions is satisfied:
\begin{enumerate}
\item $M \subseteq R$,
\item $\CC$ is principal, or
\item $R$ is of finite type over an $F$-finite field.
\end{enumerate}
However, we expect that test modules exist whenever $R$ is $F$-finite and $M$ is finitely generated, or even more generally. The first case ($M\subseteq R$) is similar to the classical proof of the existence of test elements, and hence test ideals for $F$-finite rings. Already this case extends all cases where the existence was shown in \cite{schwede_test_2009}. The second case is an application of \cite{BliBoe.CartierFiniteness}. The final case, which is given in \autoref{s.elementary}, is a by-product of a new and completely elementary approach to test ideals/modules following an idea of Anderson \cite{AndersonL}. Before we proceed to show the existence of test modules in the first two cases, we derive some basic properties about our test modules.
\begin{proposition}
\label{t.testBasic}
Let $R$ be noetherian, $\CC$ an $R$-Cartier-algebra, and $M$ a coherent left $\CC$-module. Assume that the test module $\tau(M,\CC)$ exists.
\begin{enumerate}
\item\label{t.testBasic.a} $\tau(M,\CC)=\tau(\underline{M},\CC)$.
\item\label{t.testBasic.b} $\CC_+\tau(M,\CC)=\tau(M,\CC)$, \ie~$\tau(M,\CC)$ is $F$-pure.
\item\label{t.testBasic.c} If $N \subseteq M$ and $\CC' \subseteq \CC$, then $\underline{N}_{\CC'}  \subseteq \underline{M}_{\CC}$. Suppose that $\underline{N}_{\CC'}  = \underline{M}_{\CC}$, then $\tau(N,\CC') \subseteq \tau(M,\CC)$.\footnote{In an earlier version of this paper it was claimed that the inclusion of test modules $\tau(N,\CC') \subseteq \tau(M,\CC)$ holds unconditionally. As pointed out by Daisuke Hirose, this is not the case. The problem is that minimal primes of $\underline{N}_{\CC'}$ need not be minimal primes of $\underline{M}_{\CC}$. In light of this one can state that the inclusion for test ideals holds if the minimal primes of $\underline{N}_{\CC'}$ are contained in the minimal primes of $\underline{M}_{\CC}$, see also \autoref{r.DefTestIdealAssociated}.}
\item\label{t.testBasic.d} $\tau(M,\CC)=\tau(\underline{M}_{\CC},\CC \tensor_R R_{red})$. In particular, if the latter test module exists, then so does the former. This allows us to treat the non-reduced case by reduction to the reduced case (and for this the added generality of allowing not only $M=R$ is useful).
\item\label{t.testBasic.e} If $S$ is a multiplicative set, then $\tau(S^{-1}M,S^{-1}\CC)=S^{-1}\tau(M,\CC)$. And the former exists if the latter does. Hence, the test module commutes with localization, and can be glued together from the test modules on an open cover.
\item\label{t.testBasic.f} If $\underline{M}=M$ and if $\supp M \subseteq V(I)$ for some ideal $I$ of $R$, then $\tau(M,\CC)=\tau(M,\CC/\CC I)$. Again, the former exists if the latter does. Hence one may always reduce to the case that $\supp M$ is all of $\Spec R$.
\end{enumerate}
\end{proposition}
\begin{proof}
Part \autoref{t.testBasic.a} is immediate from the definition. For part \autoref{t.testBasic.b} we may, by \autoref{t.testBasic.a}, assume that $\CC_+M = M$. Hence, by definition of the test module, for each generic point $\eta$ of a component of $\supp M$ we have $\CC_+\tau(M,\CC)_\eta = \CC_+M_{\eta}=M_{\eta}=\tau(M,\CC)_\eta$. Hence $\CC_+\tau(M,\CC)$ is a $\CC$-submodule of $M$ which agrees with $M = \underline{M}$ on each irreducible component of its support. But the test module is minimal with respect to this property, hence the inclusion $\CC_+\tau(M,\CC) \subseteq \tau(M,\CC)$ has to be equality as claimed.

For \autoref{t.testBasic.c} first observe that $\underline{N}_{\CC'} \subseteq \underline{N}_{\CC} \subseteq \underline{M}_{\CC}$, since $\CC' \subseteq \CC$, and $N \subseteq M$. If there is an equality $\underline{N}_{\CC'}  = \underline{M}_{\CC}$ then the claimed inclusion of test ideals follows immediately from the fact that there are more $\CC'$-submodules of $\underline{N}_{\CC'}  = \underline{M}_{\CC}$ than there are $\CC$-submodules.

To show \autoref{t.testBasic.d} we may, by \autoref{t.testBasic.a}, replace $M$ by $\underline{M}$ and hence assume that $\CC_+M=M$. Then the nil-radical of $R$ acts as zero, by \autoref{t.SupportReduced}. Hence we may view $M$ as a module over $R_{red}$, and replace the algebra $\CC$ by $\CC_{red}=\CC \tensor_R R_{red}$, \cf~\autoref{t.restrictionClosed} and the surrounding discussion. This shows \autoref{t.testBasic.d}.

For part \autoref{t.testBasic.e}, we again use \autoref{t.testBasic.a} and assume that $\underline{M}=M$. Since every minimal prime of $S^{-1}M$ is a minimal prime of of $M$, the fact that $\tau(M,\CC)$ agrees with $M$ for each minimal prime of $M$ by definition implies that $S^{-1}\tau(M,\CC)$ agrees with $S^{-1}M$ for each minimal prime of $S^{-1}M$. Now let $N \subseteq S^{-1}M$ be a $S^{-1}\CC$ submodule that agrees with $S^{-1}M$ for each minimal prime of $S^{-1}M$. Let $\phi: M \to S^{-1}M$ denote the localization map and let $\eta$ be a minimal prime of $M$. If $S^{-1} \subseteq  R \setminus \eta$ we have $(S^{-1}M)_\eta=M_\eta$. Since in this case $\eta$ is also a minimal prime of $S^{-1}M$ it follows from $N_{\eta}=(S^{-1}M)_\eta=M_\eta$ and $S^{-1}\phi^{-1}(N)=N$ that $\phi^{-1}(N)_{\eta}=M_\eta$. Otherwise, if $\eta \cap S^{-1}$ is not empty then the entire $\eta$-primary component of $M$ is contained in $\ker \phi$. Hence $(\ker \phi)_\eta = M_\eta$ so in particular $(\phi^{-1}(N))_\eta = M_\eta$. Now, by the minimality of $\tau(M,\CC)$ with this property it follows that $\tau(M,\CC) \subseteq \phi^{-1}(N)$ which implies that $S^{-1}\tau(M,\CC) \subseteq N$. This shows that $\tau(S^{-1}M,S^{-1}\CC)$ exists and is equal to $S^{-1}\tau(M,\CC)$.

Part \autoref{t.testBasic.f} may be viewed as a generalization of \autoref{t.testBasic.d} and follows along the same lines.
\end{proof}

\begin{remark}\label{r.DefTestIdealAssociated}
For a perfect field $k$ let $R=k[x]$ and let $\CC=R\langle C \rangle$ a principal Cartier algebra. Let $C$ act on $R$ by $x^n \mapsto x^{n/p}$ where we interpret $x^{n/p}$ as zero whenever $n$ is not divisible by $p$ (note that this is $\kappa \cdot x^{p-1}$ where $\kappa$ denotes the classical Cartier operator, \ie~a generator of $\Hom(F_*R,R)$). Then it is easy to check that $(x) \subseteq k[x]$ is the only non-trivial $\CC$-submodule of $k[x]$, hence $\tau(k[x],\CC)=(x)$. Consider $k\cong k[x]/(x)$ with the $\CC$-module structure induced from the one on $\CC$. One immediately verifies that $\tau(k,\CC)=k$. Consider the direct sum $k[x] \oplus k$, then $\tau(k[x] \oplus k,\CC)=(x) \oplus 0$ since this is the smallest $\CC$-submodule of $k[x] \oplus k$ that for the only minimal prime $(0)$ of $k[x] \oplus k$ agrees with $k[x] \oplus k$ after localization at $(0)$. As one can see from this example, the formation of test modules does not commute with direct sums, nor does it preserve inclusions in general (take the inclusion $0 \oplus k \subseteq k[x] \oplus k$). This is due to the fact that one ignores the non-minimal associated prime ideals in the construction.

At this point one might be tempted to alter the definition of the test module of a pair $(M,\CC)$, say by defining $\tau'(M,\CC)$ to be the smallest $\CC$-submodule $N$ of $M$ such that for each associated prime $\eta$ we have $N_{\eta}=(\underline{M}_{\CC})_\eta$. This might be a valid alternative and a similar theory might develop using this definition. However, continuing the example above we get $\tau'(k[x],\CC)=(x)$, $\tau'(k,\CC)=k$ but $\tau'(k[x]\oplus k,\CC)=k[x] \oplus k$. Even though this alternative definition will preserve inclusions, it still does not preserve direct sums.

These observations suggest two valid alternatives for the definition of the test module which may have better functorial properties: Firstly, we may simply require that the test module is the smallest submodule such that $\tau(M,\CC)=\underline{M}_{\CC}$ for each generic point of $\Spec R$ itself. Secondly, as suggested to me by Karl Schwede, we may alternatively define $\tau(M,\CC)$ as the smallest submodule $N$ such that for each \emph{associated} prime $\eta$ of $\underline{M}_{\CC}$ we have $H^0_\eta(N_\eta) = H^0_\eta((\underline{M}_{\CC})_\eta)$. In both cases a similar theory to the one worked out in this paper should develop.
\end{remark}

\subsection{$F$-purity and $F$-regularity}
Above we introduced the notion of $F$-purity for pairs $(M,\CC)$. In this section we will relate it to the classical notion where $M=R$ is a reduced ring, and also introduce the concept of (strong)\footnote{In accordance with established terminology, the adjective \emph{strong} should not be omitted, but for the sake of brevity I will do so anyways.} $F$-regularity in our context. Furthermore, we will show, that under the assumption that test modules exist for all pairs $(M,\CC)$ ($\CC$ is fixed and $M$ varies) then all coherent $F$-pure $\CC$-modules have finite length in the category of $F$-pure modules. This is an extension of one of the main results in \cite{BliBoe.CartierFiniteness}.
\begin{definition}
\label{d.Fpure}
The pair $(M,\CC)$ is called \emph{$F$-regular} if and only if $(M,\CC)$ is $F$-pure and there is no proper non-zero $\CC$-submodule $N \subseteq M$ that generically agrees with $M$, \ie~for each generic point $\eta$ of a component of $\supp M$ one has $N_\eta=M_\eta$.
\end{definition}
From the definition of the test module we obtain:
\begin{lemma}
The pair $(M,\CC)$ is $F$-regular if and only if $M = \tau(M,\CC)$.
\end{lemma}
\begin{proof}
The condition that $M=\tau(M,\CC)$ implies that $\CC_+M=\CC_+\tau(M,\CC)=\tau(M,\CC)=M$ by \autoref{t.testBasic} \autoref{t.testBasic.a}. Hence $(M,\CC)$ is $F$-pure. By definition, the test module is the smallest $\CC$-submodule of $M$ that generically agrees with $\underline{M}=M$, the equality being due to $F$-purity. Hence if the smallest, namely $\tau(M,\CC)$ is equal to $M$, then there are no $\CC$ submodules of $M$ that generically agree with $M$ besides $M$ itself. Hence $(M,\CC)$ is $F$-regular.

Conversely, if $(M,C)$ is $F$-regular, then $\underline{M}=M$ by $F$-purity (which is part of its definition), and the smallest $\CC$-submodule $N$ agreeing with $M$ after localizing at each minimal prime of $M$ is $M$ itself; hence $\tau(M,\CC)=M$ and, in particular, the test module exists.
\end{proof}
Tautological consequences of these definitions are that for any pair $(M,\CC)$ the pair $(\underline{M},\CC)$ is $F$-pure, and the pair $(\tau(M,\CC),\CC)$ is $F$-regular. Furthermore, we see that, provided the test ideal exists, the collection of all $F$-pure submodules of $M$ that generically agree with $\underline{M}$ has a maximal element ($\underline{M}$ itself) and a minimal element, the test ideal $\tau(M,\CC)$. In \autoref{s.elementary} we will show that if $R$ is of finite type over an $F$-finite field, all chains of $F$-pure submodules of $M$ are finite.

The following two propositions show the relation to the definition of these notions as given by Schwede in \cite{schwede_test_2009}.
\begin{proposition}
\label{t.FpureCharacterization}
The pair $(R,\CC)$ is $F$-pure iff there is $\phi_e \in \CC_e$ for $e \geq 1$ with $\phi_e(R)=R$ (this latter condition is the definition of \emph{sharply $F$-pure} for a pair $(R,\CC)$ in \cite[Definition 3.11]{schwede_test_2009}).
\end{proposition}
\begin{proof}
    Clearly, the existence of such a $\phi_e \in \CC_e$ implies that $\CC_eR=R$ for some $e \geq 1$, hence $\CC_+R=R$ and the $F$-purity of the pair $(R,\CC)$ follows.
    For the converse, since $\CC_+R=\sum_e \sum_{\phi \in \CC_e} \phi(R)$ is an ideal, it is finitely generated, hence $\CC_+R=\sum_{i=1}^k \phi_i(R)$ where $\phi_i \in \CC_{e_i}$ for $e_i \geq 1$. So, if $\CC_+R=R$, then
    \[
        \CC_+R=\sum_{i=1}^k \phi_i(R)=R.
    \]
    Hence for each point $x \in \Spec R$ we have an index $i$ such that $\phi_i(R_x)=R_x$ since a sum of proper ideals in a local ring is always a proper ideal. But then $\phi_i^n(R_x)=R_x$ for all $n \geq 0$. We claim that $\sum_{i=1}^k\phi^{e/e_i}_i(R)=R$ where $e=e_1\cdot\ldots\cdot e_k$. This equality can be checked locally: for each point $x \in \Spec R$ we have an index $i$ such that $\phi_i^{e/e_i}(R_x)=R_x$. Now, we set $\psi_i=\phi_i^{e/e_i}$ and note that each $\psi_i \in \CC_e$. There exist $r_i \in R$ such that $\sum_{i=1}^k \psi_i(r_i)=1$. But we may view $\phi_e := \sum_{i=1}^k \psi_i\cdot r_i$ itself as an element of $\CC_e$ and then $\phi_e(1)=\sum_{i=1}^k \psi_i(r_i)=1$. Hence $\phi_e(R)=R$.
\end{proof}

\begin{proposition}
The pair $(R,\CC)$ is $F$-regular iff for every $d \in R^\circ$ there is $\phi_e \in \CC_e$ for some $e \geq 1$ with $\phi_e(dR)=R$ (this latter condition is the definition of \emph{strongly $F$-regular} for a pair $(R,\CC)$ in \cite[Definition 3.12]{schwede_test_2009}).
\end{proposition}
\begin{proof}
Suppose that for every $d \in R^\circ$ there is $\phi_e \in \CC_e$ for some $e \geq 1$ with $\phi_e(dR)=R$. In particular, for every $d \in R^\circ$ we have $\CC_+(dR)=R$. For $d=1$, this shows that $(R,\CC)$ is $F$-pure, hence, by definition $\tau(R,\CC)$ is equal to the smallest $\CC$-submodule of $R$ which agrees with $R$ generically (\ie~for each minimal prime $\eta$ of $R$). But this just means that $\tau(R,\CC) \cap R^\circ$ is non-zero. Hence for this $d \in \tau(R,\CC) \cap R^\circ$ we have $R=\CC_+(dR) \subseteq \tau(R,\CC)$, hence $\tau(R,\CC)=R$, so $(R,\CC)$ is $F$-regular.

Conversely, if $(R,\CC)$ is $F$-regular, then we already observed that $(R,\CC)$ is $F$-pure. Let $d \in R^\circ$, hence $\CC d \supseteq (d)$ generically agrees with $R$. Hence $\tau(R,\CC) \subseteq \CC d$, by the minimality of $\tau(R,\CC)$ with this property. The $F$-regularity, which says that $\tau(R,\CC)=R$ now implies that $\CC d=R$. Applying $\CC_+$ and using $F$-purity we get $\CC_+(\CC d)=\CC_+ d = \CC_+ R = R$. From the equality $\CC_+(Rd)=R$ one proceeds similarly as in the proof of the preceding proposition to construct a $\phi \in \CC_e$ for some $e \geq 1$ with $\phi(d)=1$.
\end{proof}
\begin{remark}
In the case that $R$ is a domain, one immediately sees that $(R,\CC)$ is $F$-regular if and only if $(R,\CC)$ is $F$-pure and $R$ is a simple $\CC$-module. Note that assuming the simplicity of $R$ as a $\CC$-module, $F$-purity is equivalent to $\CC_+R \neq 0$. This means that the positively graded part of $\CC$ acts trivially on $R$. Hence $R$ has to be simple as an $R$ module, \ie~$R$ is a field. Hence we have, if $R$ is a domain but not a field, that the pair $(R,\CC)$ is $F$-regular if and only if $R$ is a simple $\CC$-module.
\end{remark}

\subsection{Existence of test modules in the $F$-finite case}
\label{s.existenceTest}
We now return to the general question of existence of test modules. Here, we first point out that test modules exist if and only if in the category of $F$-pure $\CC$-modules all objects have finite length, \cf~\cite[Theorem 4.6]{BliBoe.CartierFiniteness}. This characterization also holds for not necessarily $F$-finite rings. After this result is established we prove our existence result for test ideals, which, however uses the condition of $F$-finiteness.
\begin{theorem}
\label{t.finitelength}
Let $\CC$ be an $R$-Cartier algebra and $M$ a coherent $\CC$-module. Then $M$ has finite length in the category of $F$-pure $\CC$-modules if and only if  for all $\CC$-modules $N \subseteq M$ the test module $\tau(N,\CC)$ exists.\footnote{Below we will see that this condition is equivalent to that every coherent $F$-pure $\CC$-module $N$ is generically $F$-regular, \cf~\autoref{t.ExistenceTestModule}.}
\end{theorem}
\begin{proof}
The ascending chain condition is clear since already the underlying $R$-module has finite length. Since an $F$-pure module has reduced support by \autoref{t.SupportReduced} we may assume that $R$ is reduced. We show by induction on the dimension of the support of an $F$-pure $\CC$-module $M$ that it satisfies the descending chain condition on its $F$-pure submodules. The case of zero dimensional support is clear. Let $\underline{M} = M_0 \supseteq M_1 \supseteq M_2 \ldots$ be a descending chain of $F$-pure submodules of $M$. The descending chain of their (reduced!) supports must stabilize eventually. Hence we may assume, after truncating, that all $M_i$ have the same support. Similarly, by considering the generic ranks on each irreducible component of the support we may again, after truncating, assume that these are also constant. Hence, we may assume that for all $i$ the inclusion $M_i \subseteq M$ is generically (\ie~on the generic point of each irreducible component of the support of $M$) an equality. By the definition of the test ideal, it follows that $\tau(M,\CC) \subseteq M_i$ for all $i$. Now, the stabilization of the chain $(M_i)_{i>0}$ in $M$ is equivalent to the stabilization of the chain $(M_i/\tau(M,\CC))_{i>0}$ in $M/\tau(M,\CC)$. But the latter has strictly smaller dimensional support than $M$ by definition of the test ideal. So by induction, we are done.

Conversely the finite length of $M$ implies the finite length of all $\CC$-submodules and this easily implies the existence of test modules by considering the necessarily finite descending chain of $F$-pure submodules that generically agree with $\underline{M}$.
\end{proof}
From this finite length result it is easy to derive a theory of decomposition series for coherent $F$-pure $\CC$-modules, \cf~\cite{BliBoe.CartierFiniteness}, in particular, there is a well defined notion of length.

We give now an extension of \cite{BliBoe.CartierFiniteness} on the finiteness of annihilator ideals of quotients of a $F$-pure $\CC$-module. This is again under the assumption that test modules exist, \cf~\cite{BliBoe.CartierFiniteness,EnescuHochster.FrobLocCohom}.

\begin{theorem}
\label{t.finiteAnnihilators}
Let $\CC$ be an $R$-Cartier algebra and $M$ a coherent $\CC$-module. Then the set
\[
    \CS \defeq \{ \Ann_R \underline{M}/N | N \subseteq M \text{ a $\CC$-submodule of }\underline{M} \}
\]
is a set of radical ideals, closed under primary decomposition and finite intersection.

If moreover the test module exists for all pairs $(\underline{M}/N,\CC)$, then $\CS$ is a finite set of radical ideals consisting precisely of the intersections of the finitely many prime ideals it contains.
\end{theorem}
\begin{proof}
Replacing $M$ by $\underline{M}$ we may assume that $M$ is $F$-pure. Then for all $N \subseteq M$ a $\CC$-submodule $M/N$ is also $F$-pure, hence $\Ann_R M/N$ is radical by \autoref{t.SupportReduced}. Since generally $\Ann_R M/(N \cap N') = \Ann_R M/N \cap \Ann_R M/N'$ it is clear that $\CS$ is closed under finite intersections. To see that $\CS$ is closed under primary decomposition let $I=\Ann_R M/N$ for some $N \subseteq M$ a $\CC$-submodule. Without loss we may replace $R$ by $R/I$ and $M$ by $M/N$ and hence have to show that all minimal primes of $R$ are in $\CS$. Let $\frp$ be a minimal prime of $R$, and denote by $\frp'$ the intersection of the remaining ones. The reducedness of $R$ and the definition of $\frp,\frp'$ implies that $\frp \cap \frp'=0$ and $\Ann_R \frp' = \frp$. Take $M[\frp']=\Ann_M \frp'$ which is a $\CC$-submodule of $M$ by \autoref{t.restrictionClosed} and calculate $\Ann_R M/M[\frp']=\{f \in R | fM \subseteq M[\frp'] \} = \{f \in R | f\frp'M = 0 \} = \{ f \in R | f\frp'=0 \} = \Ann_R \frp'=\frp$.

By the first part of the theorem, it is enough to show that there are only finitely many prime ideals in $\CS$. Now let $N$ be such that $\frp=\Ann_R M/N$. Without loss we may replace $R$ by $R/\Ann_R M$ and hence assume that $\Ann_R M=0$. Note that $R$ is still reduced since $M$ is $F$-pure. Now, since there are only finitely many minimal primes of $R$, we may assume that $\frp$ is not one of them. But this implies that $\frp$ contains an element of $f \in R^\circ=R \setminus \{ \text{minimal primes of $R$} \}$. In particular $fM \subseteq N$ such that $M$ and $N$ agree generically (on $\Spec R_f$). Hence $\tau(M,\CC) \subseteq N$. Now we can finish the argument by induction on the dimension of the support of $M$. Since $N \supseteq \tau(M,\CC)$ we may replace $M$ by $M/\tau(M,\CC)$ which has strictly smaller support than $M$. By induction the set $\CS'=\{\Ann_R M/N | \tau(M,\CC) \subseteq N \subseteq M\} = \{ \Ann_R (M/\tau(M,\CC))/(N/\tau(M,\CC)) | N/\tau(M,\CC) \subseteq M/\tau(M,\CC)\}$ is finite.
\end{proof}

The following theorem reduces the existence of the test module to a simpler statement which will be treated below in some cases.
\begin{theorem}
\label{t.ExistenceTestModule}
Let $\CC$ be a $R$-Cartier-algebra and $M$ a coherent left $\CC$-module. Then the test module $\tau(M,\CC)$ exists if and only if on some open subset $\Spec R_c$ for $c\in R$ whose intersection with $\supp \underline{M}_{\CC}$ is dense in $\supp \underline{M}_{\CC}$  the pair $(\underline{M_c},\CC_c)$ is $F$-regular.

In this setting, the test module $\tau(M,\CC)$ is equal to the $\CC$-submodule of $M$ generated by $c^t\underline{M}$, where $t \in \mathbb{N}$ may be chosen arbitrarily.
\end{theorem}
\begin{proof}
By \autoref{t.testBasic} \autoref{t.testBasic.a} we have $\tau(\underline{M}_{\CC},\CC)=\tau(M,\CC)$ such that we may replace $M$ by $\underline{M}_{\CC}$ and assume that $(M,\CC)$ is $F$-pure. By \autoref{t.testBasic} \autoref{t.testBasic.f} we may replace $R$ by $R/\Ann_R M$ and hence assume that $\supp M=\Spec R$ and that $R$ is reduced. Now assume that there is an affine dense open subset $\Spec R_c \subseteq \Spec R$ such that the pair $(M_c,\CC_c)$ is $F$-regular, \ie~$\tau(M_c,\CC_c)=M_c$. Let $N \subseteq M$ be a $\CC$-submodule such that $N_{\eta}=M_{\eta}$ for each generic point of a component of $\Spec R$. We claim that $cM \subseteq N$. This would finish the argument since then, the $\CC$-submodule of $M$ generated by $cM$ is the test module $\tau(M,\CC)$. To see the claim we observe that the inclusion $N_c \subseteq M_c$ has to be equality due to the $F$-regularity of $M_c$. This means there is a $t \geq 0$ such that $c^tM \subseteq N$, since $M$ and hence $N$ is assumed to be finitely generated over $R$. Hence we have $c^t(M/N)=0$.  Since $N$ is a $\CC$-submodule of $M$, the quotient $M/N$ is also a $\CC$-module. Since $\CC_+M=M$ it follows that $\CC_+(M/N)=M/N$, hence $\Ann_R M/N$ is radical by \autoref{t.SupportReduced}. If follows that in fact $c(M/N)=0$, hence $cM \subseteq N$ as claimed. This finishes one implication of the first claim, and, along the way, the second claim.

For the converse we assume that the test module $\tau(M,\CC)$ exists. As above we may assume that $M$ is $F$-pure and $\supp M = \Spec R$. By definition the test module $\tau(M,\CC)$ agrees with $M$ on each component of $\Spec R$. Hence we can find a $c \in R$, not contained in any minimal prime of $R$ such that $M_c = \tau(M,\CC)_c = \tau(M_c,\CC_c)$, which shows that $M_c$ is $F$-regular.
\end{proof}
Hence, to show the existence of the test module for a pair $(M,\CC)$ it is enough to show that the pair $(\underline{M},\CC)$ is generically $F$-regular, \ie~there is an open subset $U=\Spec R_c \subseteq \Spec R$ such that $\Spec R_c \cap \supp {\underline{M}}$ is dense in $\supp \underline{M}$ and the pair $(\underline{M}_c,\CC_c)$ is $F$-regular. In fact this statement may be further reduced to a much simpler statement. First one may replace $M$ by $\underline{M}$ and hence assume that $(M,\CC)$ is $F$-pure. The $F$-purity implies that $I=\Ann_R M$ is a radical ideal, and we may replace $R$ by $R/I$, $\CC$ by $\CC/\CC I$ and hence assume that $R$ is reduced, and $\supp M = \Spec R$. Since $R$ is reduced, we may replace $R$ by any $R_c$ such that $\Spec R_c$ is a dense subset of the regular locus of $\Spec R$ (which is a dense open subset of $\Spec R$). In this case $R_c$ is the finite product of regular domains and we may consider each factor separately. Thus we may assume that $R$ is a regular domain. We can localize further and assume that $M$ is free of finite rank.

Furthermore one may also replace the Cartier algebra in some cases without changing the test ideal. Two observations of this type follow below.
\begin{lemma}
Let $M$ be a coherent $\CC$-module. There exists a \emph{finitely generated} Cartier subalgebra $\CC' \subseteq \CC$ such that $\underline{M}_{\CC}=\underline{M}_{\CC'}$. Furthermore, if the test module $\tau(M,\CC')$ exists, then so does $\tau(M,\CC)$ and one has $\tau(M,\CC') \subseteq \tau(M,\CC'') = \tau(M,\CC)$ for some \emph{finitely generated} Cartier algebra $\CC' \subseteq \CC'' \subseteq \CC$.
\end{lemma}
\begin{proof}
By replacing $M$ by $\underline{M}_{\CC}$ we may assume that $M$ is $F$-pure. Then  $\CC_+M=\sum_{e \geq 1}\sum_{\phi \in \CC_e} \phi(M) = M$ and since $M$ is finitely generated, there are finitely many $\phi_i \in \CC_{e_i}$ for $i=1\ldots k$, and $e_i \geq 1$ such that $\sum_{i=1}^k \phi_i(M)=M$. Taking $\CC'$ as the $R$-subalgebra generated by these $\phi_i$ we see that also $\CC'_+(M)=M$. Supposing that $\tau(M,\CC')$ exists, it follows by the criterion given in \autoref{t.ExistenceTestModule} that $\tau(M,\CC)$ exists as well. With $c$ as in \autoref{t.ExistenceTestModule} we have $\tau(M,\CC)=\CC cM$. But since $\tau(M,\CC)$ is finitely generated, $\CC cM=\sum_{i=1}^t \psi_i cM$ for some $\psi_i \in \CC_{e_i}$. Enlarging $\CC'$ to $\CC''$ by the finitely many $\psi_i$ we achieve that $\tau(M,\CC'')=\tau(M,\CC)$ as claimed.
\end{proof}
\begin{lemma}
Let $\CC' \subseteq \CC$ be $R$-Cartier-algebras. Let $M$ be a coherent $\CC$-module such that the pair $(M,\CC')$, and hence the pair $(M,\CC)$, is $F$-pure. Then, if the test module $\tau(M,\CC')$ exists, then so does $\tau(M,\CC)$. Furthermore, if $\CC=\sum_{i=1}^m \CC'\phi_i$ for some $\phi_i \in \CC_{e_i}$ (\ie~$\CC$ is finitely generated as a left $\CC'$-module), then $\tau(M,\CC')=\tau(M,\CC)$.
\end{lemma}
\begin{proof}
The existence of the test module is equivalent to the existence of an open subset $\Spec R_c$ such that $(M_c,\CC'_c)$ is $F$-regular and $\Spec R_c$ is dense in $\supp M$. Since the $F$-regularity of $(M_c,\CC'_c)$ implies that of $(M_c,\CC_c)$ we see that the test ideal $\tau(M_c,\CC_c)$ also exists, \cf~\autoref{t.ExistenceTestModule}. Furthermore, we have $\tau(M,\CC')=\CC' c^tM$ for all $t \geq 1$, and similarly for $\CC$. Taking $t=p^e$ for $e=e_1\cdot\ldots\cdot e_m$ we get
\[
    \tau(M,\CC)=\CC c^{p^e}M = \sum_{i=1}^m \CC'\phi_i c^{p^e}M = \sum_{i=1}^m \CC' c^{p^{e/e_i}} \phi_i M \subseteq \sum_{i=1}^m \CC' c M = \CC' cM = \tau(M,\CC').
\]
The reverse inclusion is obvious by \autoref{t.testBasic} \autoref{t.testBasic.c}.
\end{proof}

Summarizing, we have reduced the existence of test modules to the following statement:
\begin{proposition}
The existence of test modules follows once one shows that for each \emph{regular domain} $R$, and \emph{finitely generated} $R$-Cartier-algebra $\CC$, and \emph{$R$-free} left $\CC$-module $M$, if the pair $(M,\CC)$ is $F$-pure, then, after possibly localizing, it is $F$-regular.
\end{proposition}
We will show this statement in two important special cases:
\begin{proposition}
\label{t.F-pureIsGenericallyF-regular}
Assume that $R$ be an \emph{$F$-finite} regular domain. Let $\CC$ be a $R$-Cartier-algebra, and $M$ an $R$-free coherent left $\CC$-module. If
\begin{enumerate}
\item\label{t.F-pureIsGenericallyF-regular.a} either $\rank M \leq 1$, or
\item\label{t.F-pureIsGenericallyF-regular.b} $\CC$ is principally generated
\end{enumerate}
then $(M,\CC)$ is generically $F$-regular.
\end{proposition}
\begin{proof}
Since $R$ is regular and $F$-finite, $\Fr_*R$ is flat over $R$, and hence locally free of finite rank over $R$, by \cite{Kunz}. Hence after further localizing we may assume that $\Fr_*R$ is a free $R$-module of finite rank. In \autoref{t.F-pureIsGenericallyF-regular.a} we have that $M \cong R$, hence we may also replace $\CC$ by its image in $\CC_R$. Now, take any $\phi_e \in \CC_e$ nonzero for some $e \geq 1$. Since $\CC_{R,e}$ is a free rank $1$ right $R$-module (here we use $F$-finiteness), the nontrivial submodule $\phi_e R \subseteq \CC_{R,e}$ agrees with $\CC_{R,e}$ after further localizing. Hence we may assume that $\CC_{R,e} = \phi_e R$.
We replace $\CC$ by the principal subalgebra generated by $\phi_e$ since, if we show $F$-regularity for this smaller algebra, then it follows for the original one (at this point the proof of \autoref{t.F-pureIsGenericallyF-regular.a} reduces to \autoref{t.F-pureIsGenericallyF-regular.b}, but we nevertheless continue with the easy argument for \autoref{t.F-pureIsGenericallyF-regular.a}). Let $0 \neq x \in R$ be arbitrary. Since $\Fr^e_*R$ is free over $R$ there is a $\psi \in \CC_{R,e}=\CC_e$ with $\psi(x)=1$. Hence the $\CC$-submodule of $R$ generated by $x$ is all of $R$. Hence $R$ is simple as a $\CC$-module, hence, since we assumed that $(R,\CC)$ is also $F$-pure, the pair $(R,\CC)$ is $F$-regular.

Part \autoref{t.F-pureIsGenericallyF-regular.b} is an immediate consequence of \cite[Corollary 4.4]{BliBoe.CartierFiniteness}.
\end{proof}
As an immediate corollary we obtain from this
\begin{corollary}
\label{t.existenceTestModule2}
Assume that $R$ is $F$-finite. Let $\CC$ be a $R$-Cartier-algebra and $M$ a coherent left $\CC$-module. In the following cases the test module $\tau(M,\CC)$ exist.
\begin{enumerate}
\item $M \subseteq R$ or, more generally, if $\rank \underline{M}_\eta \leq 1$ for each generic point $\eta$ of an irreducible component of $\supp \underline{M}$.
\item $\CC$ is principally generated.
\end{enumerate}
\end{corollary}
\begin{proof}
Both statements follow from \autoref{t.F-pureIsGenericallyF-regular} via the reduction steps that were outlined just before \autoref{t.F-pureIsGenericallyF-regular}.
\end{proof}
Let us record how the generality of our setup provides a way to define the test ideal for not necessary reduced rings $R$: as argued before the test ideal $\tau(R,\CC)$ stays the same if we replace the pair $(R,\CC)$ by the $F$-pure pair $(\underline{R},\CC)$. Since the support of $\underline{R}$ is reduced by \autoref{t.SupportReduced} we may view the pair $(\underline{R},\CC_{red})$ as a pair over the reduced ring $R_{red}$. Now the test module $\tau(\underline{R},\CC_{red})$ is clearly a $R_{red}$-submodule of $\underline{R} \subseteq R$, hence may be viewed as an ideal of $R$. Note that for this construction our more general setup of constructing test modules for pairs $(M,\CC)$ where the first entry is any finitely generated $R$-module is very convenient. We summarize this construction in a theorem.
\begin{theorem}\label{t.test-nonreduced}
Let $R$ be a not necessaryly reduced, noetherian ring and $\CC$ a Cartier algebra on $R$. Then the test ideal $\tau(R,\CC)$ is equal to the test module of the pair $(\underline{R},\CC_{red})$ over the reduced quotient $R_{red}$ via the natural inclusions $\tau(\underline{R},\CC_{red}) \subseteq \underline{R} \subseteq R$.
\end{theorem}

\begin{example}
To illustrate our definition of the test ideal in the non-reduced case, let us compute a simple example. Consider $R=\frac{k[x,y]}{(x^2y)}$ where $k$ is a perfect field. This ring is not reduced and its reduced quotient is equal to $R_{red}=\frac{k[x,y]}{(xy)}$. The Cartier operator on $R$ is given by sending a monomial $x^iy^j$ to $x^{\frac{i-1}{p}+1}y^{\frac{j}{p}}$ where a non-integer in an exponent renders the expression equal to zero. From this one checks easily that $\underline{R}=\frac{(x)}{(x^2y)}$. As an $R_{red}$-module this is isomorphic to $R_{red}$ itself by sending $1 \in R_{red}$ to $x \in \frac{(x)}{(x^2y)}$. Now, it is easy to check that the test ideal of $R_{red}$ is equal to the maximal ideal $(x,y)$. Hence, the test ideal of $R$ is just the ideal $(x^2,xy)$.
\end{example}

\begin{remark}
In a similar fashion as we are treating the non-reduced case one can avoid the very reasonable condition of non-degeneracy which is used in \cite{schwede_test_2009}. We accomplish this again by our key observation that the test module $\tau(M,\CC)$ is equal to the test module $\tau(\underline{M},\CC)$ of the $F$-pure pair $(\underline{M},\CC)$.
\end{remark}

\forget{%
\begin{remark}
Let us remark on the relation of our setup with a non-degeneracy condition that Schwede introduced in \cite{schwede_test_2009}.
\begin{definition}[\cite{schwede_test_2009}]
Let $R$ be a noetherian ring.
\begin{enumerate}
        \item Let $\phi_e: M \to N$ be $p^{-e}$-linear. The \emph{support} of $\phi_e$ is the set of points $x \in \Spec R$ such that the stalk $(\phi_e)_x$ is non-zero. Equivalently, this is the support of the $R$-submodule $\phi_e(M)\subseteq N$. If $N$ is coherent, the support is closed.
        \item Let $\phi_e: R \to R$ be $p^{-e}$-linear. The \emph{locus of degeneracy} of $\phi_e$, denoted by $\degloc(\phi_e)$, is the set of all $x \in X$ such that $(\phi_e)_x \cdot R_x \neq \End_e(R_x)$, \ie the locus where $\phi_e$ is \emph{not} a (right) $R$-generator of $\End_e(R)$. We denote the complement of $\degloc(\phi_e)$ by $\ndsupp(\phi_e)$ and call it the \emph{non degenerate support} of $\phi_e$. If $R$ is coherent and $X$ is $F$-finite, then $\degloc(\phi_e)$ is a closed, and $\ndsupp(\phi_e)$ is an open set.\m?{What is the corresponding notion for general $M$?}
\end{enumerate}
\end{definition}
From the claims that are made in this definition we only point out that the closedness of the degeneracy locus follows once one notices that $\End_e(R)=\Hom_R(\Fr^e_*R,R)$ is coherent $\Fr$ is a finite map.
\begin{definition}
    Let $\phi_e \in \CC_{R,e}$. Then $\phi_e$ is called \emph{non-degenerate} if one of the following equivalent conditions is satisfied.
    \begin{enumerate}
        \item\label{z.d.nd.itemA} $\phi_e(R) \cap R^\circ \neq \emptyset$.
        \item\label{z.d.nd.itemB} $\phi_e$ is generically non-zero, \ie~its support is dense, or equivalently $(\phi_e)_\eta \neq 0$ for every minimal prime of $R$.
        \item\label{z.d.nd.itemD} $\phi_e$ is generically non-degenerate, \ie the non-degenerate support of $\phi_e$ is dense, or equivalently, $(\phi_e)_\eta \cdot R_\eta = \End_e(R_\eta)$.
        \item\label{z.d.nd.itemC} The map $R\langle \Ca^e\rangle \to \CC_R$ sending $\Ca^e$ to $\phi_e$ is generically injective.
        \item\label{z.d.nd.itemE} The map $R\langle \Ca^e\rangle \to (\CC_R)_{e*}=\dirsum_{n\geq 0}\CC_{en}(R)$ sending $\Ca^e$ to $\phi_e$ is generically an isomorphism.
    \end{enumerate}
\end{definition}
\begin{proof}
Assuming \autoref{z.d.nd.itemA} there is $r \in R$ such that $\phi_e(r)=s \not\in R^\circ$. Then on $U =\Spec R_s$ we have $\phi_e(r/s^q)=1$, such that $\phi_e(R_s)=R_s$. So in particular, the support of $\phi_e$ contains the dense set $U$ and we have shown \autoref{z.d.nd.itemB}. \\
Assuming \autoref{z.d.nd.itemD} take $U$ such that on $\phi_e \cdot R=\End_e(R)$ on $U$. Possibly shrinking $U=\Spec R_c$ for some $c \not\in R^\circ$ such that $R_c$ is regular it follows that the evaluation map $\End_e(R_c) \to R_c$ is surjective. But this implies that $\phi_e(R_c)=R_c$ and hence the $\phi_e(R)$ contains a power of $c \not\in R^\circ$, and \autoref{z.d.nd.itemA} follows. \\
Since $\supp(\phi_e) \supseteq \ndsupp(\phi_e)$ one has the implication \autoref{z.d.nd.itemD}$\Rightarrow$\autoref{z.d.nd.itemB}. Clearly \autoref{z.d.nd.itemE}$\Rightarrow$\autoref{z.d.nd.itemC} and \autoref{z.d.nd.itemC}$\Rightarrow$\autoref{z.d.nd.itemB}. \\
We start with \autoref{z.d.nd.itemD} $\Rightarrow$ \autoref{z.d.nd.itemE}: Let $U=\ndsupp(\phi_e)$ which we assume to be dense. That means on $U$ we have $\phi_e \cdot R = \End_e(R)$. Since $\End_e(R)=\Hom_{R}(\Fr^e_*R,R)$ is torsion free\m?{Being the dual of a coherent. Or maybe just shrink $U$ so that it is Gorenstein such that $\End_e(\CO_U)$ is even locally free} this shows that $\End_e(R)$ is actually free $R \cong \phi_e \cdot R = \End_e(R)$. This implies that $\End_{ne}(R)$ is also free, with generator $\phi^n_e$.\footnote{This is shown in \cite{schwede_f-adjunction_2009}. Let $R \to S$ be finite. $M$ coherent $S$ module. Assume that $\Hom_R(S,R) \cong S$, then the natural map
\[
    \Hom_R(S,R) \tensor_S \Hom_S(M,S) \to \Hom_R(M,R)
\]
given by composition is a isomorphism (of right $S$-modules, as one should always think of $\Hom_R(S,R)$ as an $R-S$-bimodule, and this is how we understand the tensor). This comes down to using adjointness of $\Hom$ and $\tensor$ and some keeping track of the maps. Firstly
\[
    \Hom_R(S,R) \tensor_S \Hom_S(M,S) \cong S \tensor_S \Hom_S(M,S) \cong \Hom_S(M,S)
\]
sending $\alpha \tensor \psi$ to $\psi$, where $\alpha$ is the image of $1$ under the assumed isomorphism $S \cong \Hom_R(S,R)$. Now,
\[
    \Hom_S(M,S) \cong \Hom_S(M,\Hom_R(S,R)) \cong \Hom_R(M\tensor_S S,R) \cong \Hom_R(M,R)
\]
where the map is given by sending $\psi \mapsto \alpha \cdot \psi \mapsto \alpha \circ \psi$. This implies that if $\Hom_S(M,S)$ is $S$-generated by some maps $\psi_i$, then $\Hom_R(M,R)$ is generated by $\alpha \circ \psi_i$. Applying this to $S=\Fr^e_*R$ and $M=\Fr^{ne}_*R$ our claim follows} Hence for all $n$ we have (on $U$) isomorphisms $R \cong \phi_e^n \cdot R \cong \End_e(R)$, which is the statement of $\autoref{z.d.nd.itemE}$. \\
It remains to show the implication \autoref{z.d.nd.itemB} implies \autoref{z.d.nd.itemD}: Let $U$ be the support of $\phi_e$. Further shrinking $U=\Spec R_c$ (say, to be contained in the Gorenstein locus of $X$) we have that $\End_{e}(R_c)$ is locally free of rank 1. But inclusion $\phi_e \cdot R_c \subseteq \End_*(R_c)$ of rank $1$ torsion free $R_c$-modules is generically an isomorphism, hence $\phi_e$ is generically non-degenerate.
\end{proof}
If $R$ is a domain, then the only degenerate $\phi_e \in \End_e(R)$ is the zero map. So in this case the degenerate case is the trivial case.
\end{remark}
} 

\subsection{Cartier algebras arising from pairs and triples}
\label{s.CartierAlgebrasExamples}
We quickly recall from \cite{schwede_test_2009} how the pairs and triples in birational geometry fit into our context. We choose a framework that is suitable to work in positive characteristic and define:
\begin{definition}
An \emph{$F$-graded system} of ideals $\fra^\centerdot=(\fra^a)_{a\in \mathbb{N}}$ is a descending sequence indexed by the positive integers such $(\fra^a)^{[p^b]} \cdot \fra^b \subseteq \fra^{a+b}$.
\end{definition}
Whenever one has a \emph{graded system} of ideals $\fra_\centerdot$, \ie~a descending sequence satisfying $\fra_a \cdot \fra_b \subseteq \fra_{a+b}$ (\cf~\cite[Definition 1.1]{ein_uniform_2001}) one obtains an \emph{$F$-graded system} simply by setting $\fra^b:=\fra_{p^b-1}$. This is easy to check. A key example of ($F$-)graded systems are those that arise from real powers of an ideal, \ie~$\fra^a=I^{\lceil t(p^a-1) \rceil}$, where $I$ is some ideal and $t$ a non-negative real number. By abuse of language we denote this $F$-graded system simply by $I^t$.

\begin{definition}
Let $\CC$ be an $R$-Cartier-algebra. If $\fra^\centerdot$ is a $F$-graded system, then we define the $R$-Cartier-algebra $\CC^{\fra^\centerdot}$ to be given by $(\CC^{\fra^\centerdot})_e \defeq \CC_e \cdot \fra^{e}$.
\end{definition}
It is clear that this defines a $R$-Cartier-algebra.

If $R$ is normal, let $\Delta$ be an effective $\mathbb{Q}$-divisor. In the case that $\CC \subseteq \CC_R$ Schwede associates in \cite[Remark 3.10]{schwede_test_2009} to the pair $(R,\Delta)$ the Cartier algebra $\CC^\Delta \subseteq \CC$ defined as follows: Given $\phi \in \CC^e$ this is a map $F^e_*R \to R$. Tensoring with $R(\lceil \Delta \rceil)$ we get an induced map
\[
    F^e_*R \tensor R(\lceil \Delta \rceil) \to[\phi \tensor \id] R(\lceil \Delta \rceil)\, .
\]
The reflexive hull of the left hand side is by the projection formula isomorphic to $F^e_*R(p^e\lceil \Delta \rceil)$. Since $R(\lceil \Delta \rceil)$ is reflexive the map $\phi \tensor \id$ hence factors through $F^e_*R(p^e\lceil \Delta \rceil)$. Composing the resulting map $F^e_*R(p^e\lceil \Delta \rceil) \to R(\lceil \Delta \rceil)$ with the inclusion $F^e_*R(\lceil (p^e-1)\Delta \rceil) \subseteq F^e_*R(p^e\lceil \Delta \rceil)$ we obtain a map $\phi': F^e_*R(\lceil (p^e-1)\Delta \rceil) \to R(\lceil \Delta \rceil)$. Now we define that $\phi \in \CC_e^{\Delta}$ if the image of the induced map $\phi'$ lies in $R \subseteq R(\lceil \Delta \rceil)$. In other words, $\CC_e^{\Delta} = \CC_e \cap \Image(\Hom_R(F^e_*R(\lceil(p^e-1)\Delta\rceil),R) \to \Hom(F^e_*R,R))$. In this case he shows the following proposition stating the connection to the more classical notions of test ideals.

\begin{proposition}[\cite{schwede_test_2009}]
Let $R$ be $F$-finite. For a reduced ring $R$ and graded system of ideals $\fra^\centerdot$ one has $\tau_{b}(R,\fra^\centerdot)=\tau(R,\CC_R^{\fra^\centerdot})$. If $R$ is normal and $\Delta$ is an effective $\mathbb{Q}$-divisor, then $\tau_b(R,\Delta,\fra^\centerdot)=\tau(R,\CC_R^{\Delta},\fra^\centerdot)=\tau(R,(\CC_R^{\Delta})^{\fra^\centerdot})$.
In both statements, the test ideals for pairs/triples $\tau_b$ are the classically defined big test ideals.
\end{proposition}

This shows that one can reduce questions about pairs and triples $(R,\Delta,\fra^\centerdot)$ to that of pairs of the type $(R,\CC)$. However, for some considerations, such as the discreteness of jumping numbers which we discuss below, the framework of triples $(R,\CC,\fra^\centerdot)$ (which would be equivalent to the pair $(R,\CC^{\fra^\centerdot})$) is more natural to phrase them in, so we will use both notations interchangeably.

One may generalize this setup even further, and consider a tuple of graded system of ideals $\fra_1^\centerdot,\ldots,\fra_n^\centerdot$, a typical example of which would be to take $\frb_i$ to be some ideals and $\fra_i^\centerdot=\frb_i^{t_i}$ for some non-negative real numbers $t_i$. Then, analogously this yields yet another graded system $\fra^\centerdot$ defined by $\fra^e:=\fra_1^e\cdot\ldots\cdot\fra_n^e$. And we may consider the algebra $\CC^{\fra^\centerdot}$. Then the test ideal $\tau(R,\CC^{\fra^\centerdot})$ is equal to the mixed test ideal $\tau(R,\CC,\frb_1^{t_1}\cdot\ldots\cdot\frb_n^{t_n})$ in the special case when $\fra_i^\centerdot=\frb_i^{t_i}$.

\section{An elementary approach to test ideals}
\label{s.elementary}
In this section we pick up an idea from Anderson's paper \cite{AndersonL} \emph{An elementary approach to $L$-fuctions mod $p$} and give a completely elementary approach to test ideals, test modules, and, more generally, to $F$-pure pairs $(M,\CC)$ in the case that the base ring $R=S/I$ is a quotient of a polynomial ring $S=k[x_1,\ldots,x_n]$ over an $F$-finite field.

Roughly speaking, the idea is that given an $R$-Cartier-algebra $\CC$ (satisfying some boundedness condition) and coherent $\CC$-module $M$, then there is (some) finite dimensional $k$-subspace $M^\circ$, called a \emph{nukleus} for $(M,\CC)$, such that any $F$-pure submodule $N \subseteq M$ is generated by elements in $M^\circ$. In particular, this implies that there are no infinite chains of $F$-pure submodules in $M$. From the  finite length statement one can easily derive the existence of the test modules $\tau(M,\CC)$.

\subsection{Contracting property of $p^{-e}$-linear maps}

We recall, and slightly generalize, a basic result on the contracting property of $p^{-e}$-linear maps as it is given in \cite{AndersonL}, which is our key ingredient. There, the case where $k=\mathbb{F}_q$ is considered. A minor variation of the argument however works more generally whenever the field $k$ is $F$-finite (the case $k$ perfect is completely analogous).

Let $S=k[x_1,\ldots,x_n]$ be a polynomial ring in $n$ variables over the field $k$. We generally assume that $k$ is $F$-finite, i.e.~$[k:k^p]<\infty$.  For a multi-index $i=(i_1,\dots,i_n) \in \mathbb{N}^n$ we denote by $|i|=\max_j\{i_j\}$ its maximum norm. This induces an increasing filtration of $k$-subspaces of $S$, indexed by $-\infty \cup \mathbb{N}$ where $S_d$ is the $k$-subspace freely generated by the monomials $x^i=x_1^{i_1}\cdot\ldots\cdot x_n^{i_n}$ with $|i|\leq d$ and $S_{-\infty}=0$. Let $M$ be a finitely generated $S$-module and let $m_1,\ldots,m_k$ be generators of $M$. The just introduced filtration on $S$ together with this choice of generators of $M$ induces an increasing filtration, indexed by $-\infty \cup \mathbb{N}$ on $M$, given by
\begin{align*}
    M_{-\infty}&=0, \text{ and for } d \geq 0 \\
    M_d &= S_{d}\cdot \langle m_1,\ldots,m_k \rangle
\end{align*}
For $m \in M$ we write $\delta(m)=d$ iff $m \in M_d\setminus M_{d-1}$ and call $\delta=\delta_M$ a \emph{gauge} for $M$. One should think of the gauge $\delta$ as a substitute for a degree on $M$, and the contracting property of $p^{-e}$-linear maps on $M$ is measured in terms of the gauge $\delta$. Note also, that $\delta(m)\leq d$ if there exists $r_i\in S_d$ such that $\sum r_i m_i = m$, \ie~$m$ can be written as a $S$-linear combination of the $m_i$ such that all coefficients are in $S_d$. $S$ itself has a gauge, induced by the generator $1$. Clearly, the gauge $\delta_S$ is just the degree given by grading the variable with the maximal norm on the exponents as introduced above.\footnote{In \cite{BliSchTakZha_DisccRat} we used instead the grading by total degree. The grading coming from the maximum norm, however, appears to be more natural, since then the monomials of degree $< p^e$ form a basis of $S=k[x_1,\ldots,x_n]$ over $k[x_1^{p^e},\ldots,x_n^{p^e}]$.}
\begin{example}
Consider the ring $R=\BF_p[x]$ and the Cartier linear map $\phi$ sending $1$ to $x^t$ and $x^i$ to $0$ for $1 \leq i \leq p-1$. The map $\phi$ has the following contracting property which we will explore for general Cartier linear maps below. Note that we may write $\phi(x^n) = x^{\frac{n}{p}+t}$ where one interprets the term as zero if the exponent is non-integral. Hence, for a polynomial $f(x)$ of degree $d$ we have that the degree of $\phi(f(x))$ is less or equal to $d/p + t$. That is, up to a constant $t$, the Cartier linear map $\phi$ divides the degree by $p$.
\end{example}
The following lemma summarizes, without proof, some obvious consequences of this definition.
\begin{lemma}
\label{l.BasicsGauge}
Let $M$ be finitely generated over $S=k[x_1,\ldots,x_n]$, and $\delta=\delta_M$ a gauge corresponding to some generators $m_1,\ldots,m_k$ of $M$. Then
\begin{enumerate}
\item $\delta(m)=-\infty$ if and only if $m=0$.
\item Each $M_d$ is finite dimensional over $k$ (since $S_{d}$ is).
\item $\bigcup_d M_d = M$ (since the $m_i$ generate $M$).
\item $\delta(m+m') \leq \max\{\delta(m),\delta(m')\}$
\item $\delta(fm) \leq \delta_S(f)+\delta_M(m)$
\end{enumerate}
\end{lemma}
\begin{remark}
 If $R=S/I$ is any quotient ring of $S$, then the $S$ generator $1_R$ of $R$ induces a gauge $\delta_R$, which we shall call the standard gauge on $R=S/I$ (of course this depends on the presentation). Note that we do not assume that $I$ is a homogeneous ideal of $S$. With this viewpoint we may seemingly generalize our setup slightly and consider finitely generated $R$-modules $M$, in which case $M_d = R_d \cdot \langle m_1,\ldots,m_k \rangle = S_{d} \cdot \langle m_1,\ldots,m_k \rangle$, so that we may view $M$ as an $R$-module or an $S$-module without changing the filtration or gauge.
\end{remark}

Let $\mathcal{K}_{q}=\{ k_i \}_{i\in K_q}$ be a basis of $k$ over $k^q$ for $q=p^e$. Note that the monomials in $S_{q-1}$, \ie~the monomials $\{x^i|\, 0\leq |i|<q\} \defeq \mathcal{S}_{q-1}$ are a basis of $k[x_1,\ldots,x_n]$ over $k[x_1^q,\ldots,x_n^q]$ and $\mathcal{K}_q\mathcal{S}_{q-1}$ is a basis of $S=k[x_1,\ldots,x_n]$ over $S^q=k^q[x^q_1,\ldots,x^q_n]$.
\begin{lemma}
\label{t.boundGenerator}
Let $f \in S=k[x_1,\ldots,x_n]$ with $\delta(f) \leq d$, then writing uniquely
\[
    f=\sum_{i \in K_q, 0\leq |j| < q} r^q_{i,j} k_ix^j
\]
one has $\delta(r_{i,j}) \leq \lfloor d/q \rfloor.$
\end{lemma}
\begin{proof}
Since the $x^j$ for $0 \leq |i| < q$ are a monomial basis of $S$ over $k[x^q_1,\ldots,x^q_n]$ we see that there cannot be cancellation in the sum $\sum_{0\leq |j| < q}(\sum_i  r^q_{i,j}k_i) x^j$. Hence $\delta(f)=\max_{0 \leq |j| < q}\{ \delta(\sum_i r^q_{i,j}k_i \cdot x^j)\}$. Similarly, since the $k_i$ form a basis of $k$ over $k^q$, one easily derives that a monomial $x^{lq}$ appears in $\sum_i r^q_{i,j} k_i$ if and only if it appears in $r^q_{i,j}$ for some $i$. This implies that
\[
    \delta(\sum_i r^q_{i,j} k_i \cdot x^j) \geq \delta(\sum_i r^q_{i,j} k_i) = \max_i\{\delta(r^q_{i,j})\} = q \cdot \max_i\{\delta(r_{i,j})\}\, .
\]
Hence $d\geq\delta(f) \geq q \cdot \max_{i,0\leq |j| < q}\{\delta(r_{i,j})\}$, and since $\delta(r_{i,j})$ is an integer, the claim of the lemma follows.
\end{proof}

\begin{proposition}[\protect\cite{AndersonL}, Proposition 3]
\label{t.bound1}
Let $M$ be finitely generated over $R=S/I$, and $\delta=\delta_M$ a gauge corresponding to some generators $m_1,\ldots,m_k$ of $M$. Let $k$ be $F$-finite, i.e. $[k:k^q]<\infty$ and let $\phi:M \to M$ be a $p^{-e}$-linear map. Then there is a constant $K$ such that for all $m \in M$:
\begin{align}
    \delta(\phi(m)) &\leq \frac{\delta(m)+K}{p^e} \leq \frac{\delta(m)}{p^{e}}+\frac{K}{p^e-1} \label{e.boundcartier}
\end{align}
\end{proposition}
\begin{proof}
By definition, we may write $m = \sum_{l=1}^k f_l m_l$ with $\deg f_l \leq \delta(m)$. For each $l$ write uniquely $f_l = \sum_{i\in K_{p^e},0\leq |j| < p^e} r_{l,i,j}^{p^e}k_ix^j$ using the notation introduced in \autoref{t.boundGenerator} -- then \autoref{t.boundGenerator} shows that $\delta r_{l,i,j} \leq \lfloor \delta(m)/{p^e} \rfloor$. Writing this out
\[
    \phi(m)=\sum_{l=1}^k \sum_{i,j} r_{l,i,j}\phi(k_ix^j m_l)
\]
we consequently obtain
\[
    \delta(\phi(m)) \leq \max_{l,j,i} \{\delta(r_{l,j,i})+\delta(\phi(k_jx^i m_l))\} \leq \lfloor \frac{\delta(m)}{p^e} \rfloor +\frac{K}{p^e}
\]
taking for $K = {p^e} \cdot \max_{l,j,i} \{\delta(k_jx^i m_l\})$ which exists, since $k$ is assumed to be $F$-finite. This shows (a sharpening of) the claimed inequality.
\end{proof}
\begin{corollary}
\label{t.boundMany}
Let $\phi_i$, $i=1,\ldots,n$ be $p^{-e_i}$ linear endomorphisms of $M$. Suppose that $K_i$ is the bound that works in the inequality \autoref{e.boundcartier} of \autoref{t.bound1} for $\phi_i$. Let $K = \max_i\{K_i\}$ and $e$ the greatest common divisor of the $e_i$. Then for all $m \in M$ we have
\[
    \delta(\phi_n\circ \ldots\circ \phi_1(m)) \leq \frac{\delta(m)}{p^{e_1+\ldots+e_n}} +\frac{K}{p^e-1}\ .
\]
\end{corollary}
\begin{proof}
The first step is an inductive application of the first inequality in \autoref{t.bound1} and the rest is just a computation (with $\alpha=\sum_{i=1}^n e_i/e$ and $e= \operatorname{gcd}_i\{ e_i \}$).
\[
\begin{split}
\delta ( \phi_n\circ \ldots \circ \phi_1(m))
    &\leq \frac{\delta(m) + K_1 + K_2p^{e_1}+K_3p^{e_1+e_2}+\ldots+K_n p^{e_1+e_2+\ldots+e_{n-1}}}{p^{e_1+\ldots+e_{n-1}}} \\
    &\leq \frac{\delta(m)+K(1+p^e+\ldots+p^{(\alpha-1)e})}{p^{\alpha n}} \\
    &= \frac{\delta(m)}{p^{\alpha e}} + K \frac{p^{\alpha n}-1}{p^{\alpha n}(p^e-1)} \\
    &\leq \frac{\delta(m)}{p^{\sum_i e_i}} +\frac{K}{p^e-1}
\end{split}
\]
\end{proof}

\subsection{$F$-pure submodules and test ideals}
For the rest of this section we fix a quotient $R=S/I$ of the polynomial ring in $n$ variables, a $R$-Cartier-algebra $\CC$, and a coherent $\CC$-module $M$. We want to study the set of $F$-pure submodules of $M$, \ie~the set of $\CC$-submodules $N$ of $M$ such that $\CC_+N=N$. For this question there is no difference in whether we view $M$ as an $R$-module, or as an $S$-module, so there would be no loss in generality if we assume that $R$ is the polynomial ring.

Let $\delta=\delta_M$ be the gauge induced by a (fixed) choice $m_1,\ldots,m_k$ of some generators of $M$ as an $R$-module. We say that a $R$-submodule $N \subseteq M$ is \emph{generated in gauge $\leq d$} if $N$ has a set of generators $n_1,\ldots,n_t$ with $\delta_M(n_i) \leq d$ for $i=1,\ldots,t$. As a first consequence of the contracting property of $p^{-e}$-linear maps we note:
\begin{lemma}
\label{t.BoundImage}
If $\phi$ is a $p^{-e}$-linear endomorphism of $M$ which satisfies $\delta(\phi(m)) \leq \frac{\delta(m)}{p^e}+A$ for some $A \geq 0$, and the $R$-submodule $N\subseteq M$ is generated in gauge $\leq d$, then the $R$-submodule $\phi(N)$ is generated in gauge $\leq \frac{d}{p^e} + A+1$.
\end{lemma}
\begin{proof}
First note that each element $k_ix^j$ of the basis $\mathcal{K}_{p^e}\mathcal{S}_{p^e-1}$ of $S$ over $S^{p^e}$ given above has gauge $\leq p^e-1$. The images of the elements of $\mathcal{K}_{p^e}\mathcal{S}_{p^e-1}$ in $R=S/I$ form a generating set of $R/R^q$, each with gauge $\leq p^e-1$ (this is the gauge on $R$ induced by the generator $1_R$). We denote this generating set by $\overline{\mathcal{K}_{p^e}\mathcal{S}_{p^e-1}}$.

If $n_1,\ldots,n_t$ is a set of $R$-generators of $N$ then
\[
    \phi(N)=\phi(\sum_{k=1}^t R n_k)=\phi(\sum_{k=1}^t\sum_{b \in \overline{\mathcal{K}_{p^e}\mathcal{S}_{p^e-1}}} R^{p^e}b n_k)=\sum_{k=1}^t\sum_{b \in \overline{\mathcal{K}_{p^e}\mathcal{S}_{p^e-1}}} R \phi(b n_k),
\]
hence $\phi(N)$ is generated by the finite set $\phi(b n_k)$ where $b \in \overline{\mathcal{K}_q\mathcal{S}_q}$ and $k=1,\ldots , t$. By assumption we have $\delta(\phi(b n_k))\leq \frac{\delta(b n_k)}{p^e}+A\leq \frac{(p^e-1)+\delta(n_k)}{p^e}+A \leq \frac{d}{p^e}+A+1$.
\end{proof}

If the $R$-Cartier algebra $\CC$ is generated by a single element in degree $e$, and if $M$ is a coherent $\CC$-module then \autoref{t.boundMany} shows that there is a gauge $\delta$ on $M$ and a constant $K$ such that for all $m \in M$ we have $\delta(\phi^n(m)) \leq \frac{\delta(m)}{p^{ne}}+\frac{K}{p^e-1}$. Note that the $\phi^n$ for $n \geq 1$ form a set of right $R$-module generators of the $R$-algebra $\CC$. This universal contracting property that a principal Cartier algebra satisfies motivates to define:

\begin{definition}
A pair $(M,\CC)$, consisting of an $R$-Cartier-algebra $\CC$ and a coherent $\CC$-module $M$, is called \emph{gauge bounded} if for each/some gauge $\delta$ on $M$ there exists a set $\{\psi_i | \psi_i \in \CC_{e_i}, e_i \geq 1\}_{i \in I}$ which generates $\CC_+$ as a \emph{right} $R$-module, and a constant $K$ such that for each  $\psi_i$, one has $\delta(\psi_i(m)) \leq \frac{\delta(m)}{p^e}+\frac{K}{p-1}$.
\end{definition}

\begin{proposition}
\label{t.finitelyGenisGaugeBounded}
Any pair $(M,\CC)$ where $\CC$ is finitely generated as an $R$-algebra and $M$ is coherent is gauge bounded.
\end{proposition}
\begin{proof}
Let $\{\phi_i | \phi_i \in \CC_{e_i}, e_i \geq 1 \}_{i \in I}$ be $R$-algebra generators of $\CC$. Then the set
\[
\{ \psi | \psi \text{ finite product of some $\phi_i$'s}\}
\]
is a set of \emph{right} $R$-module generators of $\CC$. If the index set $I$ is finite, then \autoref{t.boundMany} shows that there is a constant $K$ such that for each such $\psi \in \CC_e$ one has $\delta(\psi(m)) \leq \frac{\delta(m)}{p^e} + \frac{K}{p-1}$ as claimed.
\end{proof}
\begin{corollary}
If $(M,\CC)$ is gauge bounded by a constant $K$, and $N \subseteq M$ is generated in gauge $\leq d$, then $\CC_+^eN$ is generated in gauge $\leq \frac{d}{p^e}+\frac{K}{p-1}+1$. Consequently, if $N$ is an $F$-pure submodule, \ie~$\CC_+N=N$ then $N$ is generated in gauge $\leq \frac{K}{p-1}+1$ which is independent of $N$.
\end{corollary}
\begin{proof}
The assumption that $(M,\CC)$ is gauge bounded by $K$ means that there is a right $R$-module generating set $\{\psi_i \in \CC_{e_i}\}_{i \in J}$ of $\CC_+$ such that $\delta(\psi_i(m)) \leq \frac{\delta(m)}{p^{e_i}}+\frac{K}{p-1}$ for each $i \in J$. Since $N$ and hence $\CC_+N \subseteq M$ is finitely generated, there are $\phi_1,\ldots,\phi_k$, each of which a product of some $\psi_i$ with $i \in J$ such that $\CC^e_+N=\sum_{i=1}^k \phi_i(N)$ and $\phi_i \in \CC_{\epsilon_i}$ with $\epsilon_i \geq e$. By \autoref{t.BoundImage} we get that $\phi_i(N)$ is generated in gauge $\leq \frac{d}{p^{\epsilon_i}} + \frac{K}{p-1} + 1$. Hence $\CC_+^eN$ is generated in gauge $\leq \frac{d}{p^e}+\frac{K}{p-1}+1$, showing the first statement.

If in addition $N$ is $F$-pure, \ie~$\CC_+N=N$, then $\CC_+^eN=N$ for all $e \geq 0$, and hence, by the first part, $N$ is generated in gauge $\leq \frac{d}{p^e}+\frac{K}{p-1}+1$. But this holds for all $e \geq 1$, and since the gauge is an integer valued function, the term $\frac{d}{p^e}$ does not contribute for large $e$. Hence $N$ is generated in gauge $\leq \frac{K}{p-1}+1$.
\end{proof}

\begin{remark}
In the notation of the preceding corollary, this result shows that the finite dimensional $k$-subspace $M_{\leq \lfloor \frac{K}{p-1}+1\rfloor}$ completely determines the theory of $F$-pure submodules of $M$, as any $F$-pure submodule of $M$ has generators in $M_{\leq \lfloor \frac{K}{p-1}+1\rfloor}$.

Following \cite{AndersonL} we call a finite dimensional $k$-subspace $N$ of $M$ a \emph{nukleus} for $(M,\CC)$ if for every $m \in M$ there is an index $n \geq 0$ such that $\CC^n_+(m) \subseteq N$. Similarly as in our discussion so far one checks that if $(M,\CC)$ is gauge bounded by $K$, then $(M,\CC)$ has a finite dimensional nukleus $M_{\leq \lfloor \frac{K}{p-1}\rfloor}$. Of course, the intersection of any two nuklei is also a nukleus, so there exists a unique minimal nukleus $M^\circ$ for any gauge bounded pair $(M,\CC)$. The dimension of the minimal nukleus is hence an invariant of the pair $(M,\CC)$. What is its meaning? Note that a nukleus will generally not enjoy the property that every $F$-pure submodule has generators in it. For example, if $R=k[x]$ and $\kappa$ is the usual Cartier operator $\kappa(x^n)=x^{{n+1}/p-1}$, then the minimal nukleus of $(R,R\langle \kappa \rangle)$ is $0$. However note that for $R=k[x,x^{-1}]$ the minimal nukleus is $k \cdot x^{-1}$, since $x^{-1}$ is a fixed point for the Cartier action.
\end{remark}

We quickly derive from this that in a gauge bounded pair $(M,\CC)$ there are no infinite proper chains of $F$-pure $\CC$-submodules of $M$.
\begin{proposition}
\label{t.gaugeboundfinitelength}
Assume that $(M,\CC)$ is gauge bounded by a constant $K$. Then there are no infinite proper chains of $F$-pure $\CC$-submodules of $M$.
\end{proposition}
\begin{proof}
By the preceding corollary each $F$-pure submodule $N \subseteq M$ is generated in gauge $\leq \frac{K}{p-1}+1$. Hence each such $N$ is completely determined by its intersection with $M^\circ \defeq M_{\frac{K}{p-1}+1}$, which is a finite dimensional $k$-vectorspace. The intersection $N \cap M^\circ$ is a $k$-subspace of $M^\circ$. This implies, in particular, that there are no infinite chains of $F$-pure submodules of $M$.
\end{proof}

\begin{theorem}
\label{t.TestExistsAffineKalg}
For any pair $(M,\CC)$ where $\CC$ is a $R$-Cartier-algebra, $R$ is of finite type over a $F$-finite field, and $M$ is a coherent $\CC$-module, the test module $\tau(M,\CC)$ exists.
\end{theorem}
\begin{proof}
We may first replace $M$ by $\underline{M}_\CC$ and assume that $(M,\CC)$ is $F$-pure. The equality $\CC_+M=M$ implies that there are finitely many $\phi_i \in \CC_{e_i}$, $e_i \geq 1$ such that already $\sum \phi_i(M)=M$. If $\CC'$ denotes the $R$-Cartier subalgebra of $\CC$ generated by these finitely many $\phi_i$, we also have that the pair $(M,\CC')$ is $F$-pure. Since $\CC'$ is finitely generated, the pair $(M,\CC')$ is gauge bounded, and hence there are no infinite chains of $F$-pure (with respect to $\CC'$) submodules of $M$, by \autoref{t.gaugeboundfinitelength}. This implies that the intersection of all $F$-pure submodules (with respect to $\CC'$) that generically agree with $M$ is also an $F$-pure submodule that generically agrees with $M$. By definition, this intersection is the test module $\tau(M,\CC')$, which hence exists. But now, the test module $\tau(M,\CC)$ is the $\CC$-submodule of $M$ generated by $\tau(M,\CC')$, and hence exists as well.
\end{proof}
The existence of the test modules now implies, by \autoref{t.finitelength}, the finite length of $F$-pure modules in general, without the gauge boundedness assumption.
\begin{corollary}
Let $R$ be an affine $k$ algebra over the $F$-finite field $k$, and $\CC$ an $R$-Cartier algebra. Then all objects in the category of coherent $F$-pure $\CC$-modules have finite length.
\end{corollary}

\subsection{Discreteness of jumping numbers}
As another application of the gauge viewpoint on test ideals we give now a very simple proof of the discreteness of jumping numbers of test modules. We need a preparatory proposition.
\begin{proposition}
If the pair $(M,\CC)$ is gauge bounded (with bound $K/(p-1)$), and $\fra$ is an ideal, generated in gauge $\leq d$, $t$ is a non-negative real number, then the pair $(M,\CC^{\fra^t})$ is gauge bounded with bound $K/(p-1)+1+td$.
\end{proposition}
\begin{proof}
Let $\{\psi_i\}_{i\in I}$ be the system of right $R$-module generators of $\CC_+$ which satisfies $\delta(\psi_i(m))\leq \frac{\delta(m)}{p^{e_i}}+\frac{K}{p-1}$ for all $m \in M$. Since $\fra$ is generated in gauge $\leq d$, the ideal $\fra^{\lceil t(p^e-1) \rceil}$ is generated in gauge $\leq d t p^e+1$. Using again a set of generators of $R$ over $R^{p^e}$ each of gauge $\leq (p^e-1)$ we see as before that there is a set of right $R$-module generators of $\CC^{\fra^t}$, each satisfying the bound
\[
\delta(\phi(m)) \leq \frac{\delta(m)}{p^e}+\frac{K}{p-1}+1+dt=\frac{\delta(m)}{p^e}+\frac{K+(1+dt)(p-1)}{p-1}
\] This follows since each generator $\phi \in \CC_e^{\fra^t}$ is of the form $\psi \cdot b \cdot a$ where $b$ is one of the generators of $R$ over $R^{p^e}$ and $a$ is one of the generators of $\fra^{\lceil t(p^e-1) \rceil}$, the gauge of each satisfying the respective bound. Then one computes
\[
\begin{split}
    \delta(\phi(m))&=\delta(\psi(bam))\leq\frac{\delta(bam)}{p^e}+\frac{K}{p-1}=\frac{\delta(b)+\delta(a)+\delta(m)}{p^e}+\frac{K}{p-1} \\ &\leq \frac{p^e + dtp^e +\delta(m)}{p^e}+\frac{K}{p-1}=\frac{\delta(m)}{p^e}+\frac{K}{p-1}+1+dt.
\end{split}
\]
\end{proof}
Before we proceed we show the following Propsition which will justify the definition of jumping numbers in our generalized context.
\begin{proposition}\label{t.jumps}
    Let $R$ be $F$-finite, $\CC$ an $R$-Cartier algebra, $M$ a coherent $\CC$-module, $\fra$ an ideal in $R$, and $t \geq 0$ a real number.
    Then for all $\epsilon \geq 0$ we have
    \[
        \tau(M,\CC,\fra^t)\supseteq \tau(M,\CC,\fra^{t+\epsilon})
    \]
    with equality for sufficiently small $\epsilon >0$.
\end{proposition}
\begin{proof}
The proof is similar as in the more classical case \cite[Lemma 3.23]{BliSchTakZha_DisccRat} -- at least after some reductions have been made. First we may replace $M$ by $\underline{M}_{\CC^{\fra^{t}}}$ and assume that $(M,\CC^{\fra^t})$ is $F$-pure. We may further replace $R$ by $R/\Ann_R M$ and assume that $\Ann_R M=0$ and that $R$ is reduced; both reductions are due to \autoref{t.testBasic}.

We claim that now $\fra$ is not contained in the union of the minimal primes of $R$. Supposing otherwise we find by prime avoidance a minimal prime $\eta$ of $R$ containing $\fra$. Then the localization $R_{\eta}$ is artinian, irreducible and reduced, and hence a field. Since $\fra \subseteq \eta$ the image of $\fra$ in $R_{\eta}$ is zero. Hence $(\CC_e^{\fra^t})_\eta=\CC_e(\fra^{\lceil t(p^e-1)\rceil})_\eta=0$ for $e > 0$. Since $M_\eta$ is still $F$-pure and non-zero we have a contradiction since $(\CC_+^{\fra^t})_\eta=0$.

As usual denote by $R^\circ$ the set $R$ minus the union of its minimal primes. By \autoref{t.ExistenceTestModule} we choose $c \in R^\circ$ such that the pair $(M_c,(\CC^{\fra^t})_c)$ is $F$-regular. After replacing $c$ by $ca$ for some $a \in \fra \cap R^\circ$ we may assume that $c \in \fra$. For this reason we have that $(\CC^{\fra^s}_e)[c^{-1}] = \CC_e\fra^{\lceil s(p^e-1)\rceil}[c^{-1}]=\CC_e$ for all $e > 0$. Hence it follows that \emph{for all $s>0$} the pair $(M_c,\CC^{\fra^s})_c)$ is $F$-regular, in particular for $s=t+\epsilon$ with $\epsilon \geq 0$. By the second statement of \autoref{t.ExistenceTestModule} we conclude for all $\epsilon \geq 0$ that
\[
    \tau(M,\CC^{\fra^{t+\epsilon}}) = \CC^{\fra^{t+\epsilon}}cM \subseteq \CC^{\fra^t}cM = \tau(M,\CC^{\fra^t}).
\]
On the other hand, again applying part two of \autoref{t.ExistenceTestModule}, we get the reverse inclusion
\[
\begin{split}
    \tau(M,\CC^{\fra^{t+\epsilon}}) &\supseteq \sum_{e=1}^{E}\CC_e^{\fra^{t+\epsilon}}cM = \sum_{e=1}^{E}\CC_e \fra^{\lceil (t+\epsilon)(p^e-1)\rceil} cM \\
                           &\supseteq \sum_{e=1}^{E}\CC_e \fra^{\lceil t(p^e-1)\rceil}ccM  \qquad\qquad\qquad (\text{for }\epsilon < 1/(p^E-1))
                                    \\
                                    &= \sum_{e=1}^E \CC_e^{\fra^t} c^2M  \\
                                    &= \tau(M,\CC^{\fra^t}) \qquad\qquad\qquad (\text{for some }E \gg 0)
\end{split}
\]
for some sufficiently big $E$ and hence for all $\epsilon < 1/(p^E-1)$.
\end{proof}
\begin{definition}
Let $R$ be $F$-finite, $\CC$ an $R$-Cartier algebra, $M$ a coherent $\CC$-module, $\fra$ an ideal in $R$. Then a real number $t > 0$ is called a \emph{jumping number} if  $\tau(M,\CC,\fra^{t-\epsilon}) \varsubsetneq \tau(M,\CC,\fra^t)$ for all $\epsilon > 0$.
\end{definition}
The following result will immediately imply our discreteness result for the jumping numbers of test modules.
\begin{theorem}
\label{t.GaugeBoundDiscrete}
Let $(M,\CC)$ be gauge bounded, $\fra$ an ideal, and $T$ a positive real number. Then the set of test modules $\tau(M,\CC,\fra^t)$ for $0 \leq t \leq T$ is a finite set.

Similarly and under the same assumptions, the set of $F$-pure submodules $\underline{M}_{\CC^{\fra^t}}$ of $M$ for $0 \leq t \leq T$ is a finite set.
\end{theorem}
\begin{proof}
By \autoref{t.jumps} we have $\tau(M,\CC^{\fra^s}) \subseteq \tau(M,\CC^{\fra^t})$ for $s \geq t$. So the family of these test ideals is ordered by inclusion. As we have observed before, each $\tau(M,\CC^{\fra^t})$ is determined by its intersection with $M_{K/(p-1)+1+td}$, where $K$ is the gauge bound for $\CC$, and $d$ is such that $\fra$ is generated in gauge $\leq d$. Hence, in particular, each $\tau(M,\CC^{\fra^t})$ is determined by its intersection with $M_{K/(p-1)+1+Td}$, which is a finite dimensional $k$-vectorspace. Hence the chain of test ideals must be finite.

Since the set $\underline{M}_{\CC^{\fra^t}}$ for $t > 0$ is also ordered by inclusion the exact same argument shows the claim in this case.
\end{proof}

\begin{corollary}
Let $(M,\CC)$ be gauge bounded and $\fra$ an ideal. Then the jumping numbers of $\tau(M,\CC,\fra^t)$ are a discrete subset of $\mathbb{R}_{\geq 0}$.
\end{corollary}

\begin{remark}
Note that this in particular implies that if $\CC$ is a finitely generated $R$-algebra, then the jumping numbers of $\tau(M,\CC,\fra^t)$ are discrete for all ideals $\fra$. The $\mathbb{Q}$-Gorenstein cases considered in \cite{BliSchTakZha_DisccRat} reduce to the case when $\CC$ is generated by a single element, so our results here greatly generalize the ones in \cite{BliSchTakZha_DisccRat}. However, presently I do not have an understanding of what the condition that the Cartier algebra is finitely generated, or even worse, the condition that the pair $(M,\CC)$ is gauge bounded, means. For example, in the previously mentioned example of Katzman \cite{Katzman.NonFG} of a ring where $\CC_R$ is not finitely generated it is easy to verify that the resulting pair $(R,\CC_R)$ is gauge bounded.

The finitely generated-ness of $\CC_R$ one could speculate to be related to the finite generation of the (anti-)canonical ring. Also unclear is the relation to new discreteness results in characteristic zero \cite{urbinati_discrepancies_2010}.
\end{remark}

We conclude with pointing out that a Skoda type statement also holds in this very generalized setting, with essentially the same proof as in the classical case, see \cite{HaraTakagi.TestIdeal}.

\begin{theorem}[Skoda]
Let $R$ be $F$-finite, $\CC$ an $R$-Cartier algebra, $M$ a coherent $\CC$-module, $\fra$ an ideal in $R$, and $t \geq 0$ a real number. Then
\[
    \fra \cdot \tau(M,\CC,\fra^{t-1}) \subseteq \tau(M,\CC,\fra^t)
\]
with equality whenever $t$ is greater or equal to the number of generators of $\fra$.
\end{theorem}
\begin{proof}
As in the proof of \autoref{t.jumps} we may replace $M$ by $\underline{M}_{\CC^{\fra^{t-1}}}$ and then $R$ by $R/\Ann_R M$ such that we may assume that $R$ is reduced, $M$ is $\CC^{\fra^{t-1}}-F$-pure, and $\supp M = \Spec R$. Then the proof of \autoref{t.jumps} shows that we may pick $c \in \fra \cap R^\circ$ such that the triples $(M,\CC,\fra^{t-1})$ and $(M,\CC,\fra^{t})$ are $F$-regular on $\Spec R_c$. Hence by \autoref{t.ExistenceTestModule} we have
\[
\begin{split}
    \tau(M,\CC,\fra^t) &= \sum_{e \geq 1} \CC_e \fra^{\lceil t(p^e-1)\rceil} c \\  &\supseteq \sum_{e \geq 1} \CC_e \fra^{[p^e]} \cdot \fra^{\lceil (t-1)(p^e-1)\rceil} c = \fra \sum_{e \geq 1} \CC_e \fra^{\lceil (t-1)(p^e-1)\rceil} c = \tau(M,\CC,\fra^{t-1})
\end{split}
\]
since $\fra^{[p^e]} \cdot \fra^{\lceil (t-1)(p^e-1)\rceil} \subseteq \fra^{\lceil t(p^e-1)\rceil}$ showing the claimed inclusion.

To see the equality we observe that for $t$ greater or equal to the number of generators of $\fra$ one has $\fra^{\lceil t(p^e-1) \rceil} = \fra^{[p^e]} \cdot \fra^{\lceil (t-1)(p^e-1)-1\rceil}$ and hence
\[
\begin{split}
    \tau(M,\CC,\fra^t) &= \sum_{e \geq 1} \CC_e \fra^{\lceil t(p^e-1) \rceil} c^2 = \sum_{e \geq 1} \fra \cdot \CC_e \fra^{\lceil (t-1)(p^e-1)-1\rceil} c^2 \\
    &\subseteq \fra \cdot \sum_{e \geq 1} \CC_e \fra^{\lceil (t-1)(p^e-1)\rceil} c = \fra \cdot \tau(M,\CC,\fra^{t-1})
\end{split}
\]
as claimed.
\end{proof}

\section{Questions}
\label{s.question}

\begin{question}
The central question that needs to be addressed is the existence of the test modules in complete generality. As we explained in this paper this comes down to show that any $F$-pure pair $(M,\CC)$ is generically $F$-regular, \ie there is a $c \in R$ such that $\Spec R_c \cap \supp M$ is dense in $\supp M$. I expect this to be the case whenever $R$ is $F$-finite, or for even more general noetherian rings $R$ for which the existence of the test ideal in the classical setting is known.

One possible approach to this problem would be to reduce it to the case of a principally generated Cartier algebra where we know everything due to \cite{BliBoe.CartierFiniteness}. What one would need to show is that given any coherent $F$-pure $\CC$-module $M$ there exists a principally generated Cartier sub-algebra $\CC' \subseteq \CC$ such that the pair $(M,\CC')$ is still $F$-pure, at least generically. Contrary to what was suggested in the preprint version of this paper this is not true as one can find explicit examples of a Cartier algebra $\CC$ over a non-perfect field and an $F$-pure $\CC$-module $M$ such that the pair $(M,\CC')$ is not $F$-pure for any principal $\CC'$-subalgebra of $\CC$.

Hence in order to answer the above question one needs a different line of attack.
\end{question}


\begin{question}
If $M$ and $N$ are simple $F$-pure $\CC$-modules. Is $\Hom_{\CC}(M,N)$ finite? First note, that if $M$ is not isomorphic to $N$, then the $\Hom$ is zero. So we have to check that $\End_{\CC}(M)$ is finite. The simplicity of $M$ implies that $M$ has a unique associated prime $\frp$, and that $M$ is annihilated by $\frp$, by $F$-purity. After localizing at $\frp$, and then modding out by $\frp$ we may assume that the ground ring $R$ is a field, and it would be enough to answer this question in the field case.
\end{question}

\begin{question}
Can one show the finiteness of the $F$-pure submodules for general $\CC$? By what we have shown, there are no infinite chains of such submodules. From our \cite{BliBoe.CartierFiniteness} it follows that finiteness is true if $\CC$ is principally generated. This would yield an analog of the finiteness of compatibly split subvarieties from the theory of Frobenius splittings which does not refer to a particular splitting. A positive answer to the preceding questions would be helpful with this question as well, but things might be simpler (at least in the case $M=R$).
\end{question}

\begin{question}
In our elementary proof of the discreteness of the test ideals $\tau(M,\CC,\fra^t)$ for affine $k$-algebras over an $F$-finite field we have a condition on the boundedness of the algebra $\CC$. There are two related questions in this context which are interesting.
\begin{enumerate}
\item What are neccessary conditions on $R$ such that $\CC_R$ is bounded. A candidate for such a condition would be that the anti-canonical ring is finitely generated, as this might imply that $\CC_R$ is finitely generated and hence bounded. More precisely one may ask: Are there at all examples of rings $R$ such that the pair $(R,\CC_R)$ or $(M,\CC_R)$ for some coherent $R$-module is unbounded?
\item Is the boundedness condition for the discreteness of the jumping numbers really necessary? In the case of the existence of test ideals we could resort to a finitely generated subalgebra of $\CC$, which allowed us to show the existence without a boundedness condition. Is a similar argument possible for the discreteness?
\end{enumerate}
More generally, if $R$ is not of finite type over an $F$-finite field, can one proof discreteness, at least for finitely generated $R$-Cartier algebras $\CC$?
\end{question}

\begin{question}
Are there only finitely many ideals in the family $\tau(M,\CC,\fra^{t_1}\cdot\ldots\cdot\fra^{t_m})$ for $t_i$ each smaller than some fixed bound $T$? In the case that $R$ is an affine $k$-algebra over an $F$-finite field our proof above of the single ideal case shows that along any finite increasing path (coordinate-wise partial order on the exponents) in $\mathbb{R}^n$ there are only finitely many possible test ideals. This does not imply, however, the finiteness of all of them. In characteristic zero for multiplier ideals (and in the case that $R$ is $\mathbb{Q}$-Gorenstein) one has an even stronger statement. The partition on $\mathbb{R}^n$ given by the pieces where multiplier ideals $\CJ(R,\fra^{t_1}\cdot\ldots\cdot\fra^{t_m})$ are constant is rational and locally finite. Even in the smooth case such a statement is unknown in positive characteristic.
\end{question}

\begin{question}
Fixing a $R$-Cartier algebra (\textit{e.g.}~$\CC_R$), study the category of coherent $F$-pure $\CC$-modules. Interesting objects in this category are the test modules introduced here. Under reduction modulo $p$, these correspond to the multiplier ideals in characteristic zero. What do other $F$-pure modules correspond to in characteristic zero? The first obvious candidate is to look at $\underline{R}_{\CC}$ and this was already considered in \cite{fujino_supplements_2010} under the name of \emph{maximal non-$F$-pure ideal}.
\end{question}

\begin{question}
More a suggested investigation than a question: Develop Schwede's theory of $F$-adjunction \cite{schwede_f-adjunction_2009} and centers of $F$-purity \cite{schwede_centers_2008} in the general context of this paper. With the centers of $F$-purity we essentailly do this with our notion of $F$-pure pairs above. My impression is that the language we use here is very natural in this context so that this might also lead to a simplified and generalized treatment.
\end{question}

\begin{question}
Do the ideas and results in this paper shed some light on the question whether the strong and weak test ideals are equal. As I am told, this is still only known if the punctured spectrum is $\mathbb{Q}$-Gorenstein or if the ring in question is graded.
\end{question}

\def\cprime{$'$}
\providecommand{\bysame}{\leavevmode\hbox to3em{\hrulefill}\thinspace}
\providecommand{\MR}{\relax\ifhmode\unskip\space\fi MR }
\providecommand{\MRhref}[2]{%
  \href{http://www.ams.org/mathscinet-getitem?mr=#1}{#2}
}
\providecommand{\href}[2]{#2}

\end{document}